\newif\ifsiam
\newif\ifnummat

\siamtrue

 \nummattrue

\ifsiam
 \documentclass[final,leqno]{siamltex}
\else
  \ifnummat
     \documentclass[nummat]{svjour}
  \else
     \documentclass[twocolumn,fleqn,runningheads,final]{svjour2}
     \journalname{Computing and Visualization in Science}
     \usepackage{mathptmx}
  \fi

\fi

\usepackage{curves,amsmath,amssymb,color,fancybox,graphicx,exscale,verbatim}

\definecolor{dark}{gray}{0.6}
\definecolor{light}{gray}{0.8}

\def\eref#1{(\ref{#1})}

\def\Re{\mathbb R}

\def\ip#1#2{( #1 , #2 )}
\def\bigip#1#2{\big( #1 , #2 \big)}
\def\biggip#1#2{\bigg( #1 , #2 \bigg)}
\def\Bigip#1#2{\Big( #1 , #2 \Big)}

\def\bip#1#2{\langle #1 , #2 \rangle}
\def\bigbip#1#2{\big\langle #1 , #2 \big\rangle}

\def\norm#1{|\!| #1 |\!|}
\def\bignorm#1{\big|\!\big| #1 \big|\!\big|}

\def\snorm#1{| #1 |}

\def\enorm#1{|\!|\!| #1 |\!|\!|}

\def\rebAuthor{Randolph E. Bank}

\def\mmAuthor{Maximilian S. Metti}

\ifsiam
    \def\rebShortAuthor{R.~E.~Bank}
    \def\mmShortAuthor{M.~S.~Metti}
\else
    \def\rebShortAuthor{Bank}
    \def\mmShortAuthor{Metti}
\fi

\def\rebAddress{Department of Mathematics, University of California, San Diego,
 La Jolla, California 92093-0112. Email: rbank@ucsd.edu.}

\def\mmAddress{Department of Mathematics, University of California, San Diego,
 La Jolla, California 92093-0112. Email: mmetti@ucsd.edu.}

\def\rebThanks{The work of this author was supported by the National
Science Foundation under  contract DMS-xxxxxxx.}

\def\mmThanks{The work of this author was supported by the National
Science Foundation under  contract DMS-xxxxxxx.}

\title{Generalized time integration schemes for space-time moving finite elements}

\def\shortTitle{Moving Finite Elements with TR-BDF2}

\def\myKeywords{TR-BDF2, Moving Finite Elements, Arbitrary Lagrange-Eulerian Methods,
Method of Characteristics, Convection-Dominated, Error Analysis}

\def\myAMS{65M55, 65F10}

\def\myAbstract{
In this paper, we analyze and provide numerical illustrations for a moving finite element method applied to convection-dominated, 
time-dependent partial differential equations.
We follow a method of lines approach and utilize an underlying tensor-product finite element space that 
permits the mesh to evolve continuously in time and undergo discontinuous  reconfigurations at discrete time steps.
We employ the TR-BDF2 method as the time integrator for piecewise quadratic tensor-product spaces,  and provide an
almost symmetric error estimate for the procedure.
Our numerical results validate the efficacy of these moving finite elements.
}

\begin{document}


\ifsiam
  \author{\rebAuthor%
         \thanks{\rebAddress\ \rebThanks}
          \and
          \mmAuthor%
         \thanks{\mmAddress\ \mmThanks}
         }
  \maketitle

  \begin{abstract}\myAbstract\end{abstract}
  \begin{keywords}\myKeywords\end{keywords}
  \begin{AMS}\myAMS\end{AMS}
  \pagestyle{myheadings}
  \thispagestyle{plain}
  \markboth{\rebShortAuthor\ and \mmShortAuthor }{\shortTitle}

\else


  \ifnummat
      \author{\rebAuthor%
             \thanks{\rebThanks}
             \and
             \mmAuthor%
             \thanks{\mmThanks}
              }
  \else
      \author{\rebAuthor%
             \thanks{\rebShortAuthor : \rebThanks}
             \and
             \mmAuthor%
             \thanks{\mmShortAuthor : \mmThanks}
              }
  \fi
  \institute{\rebShortAuthor : \rebAddress \\ \mmShortAuthor : \mmAddress}
  \date{Received: \today\  / Accepted: date}
  \maketitle
  \begin{abstract}\myAbstract\end{abstract}
  \begin{keywords}\myKeywords\end{keywords}
  \begin{subclass}\myAMS \end{subclass}
\fi



\section{Introduction}\label{sec1}

For parabolic equations, the method of lines is an efficient approach for computing a numerical solution by converting the partial differential equation into
a coupled system of ordinary differential equations.
This provides a great deal of flexibility in how the solution may be computed, as the time discretization then becomes independent of the spatial discretization.
For finite element methods, the spatial dimensions are discretized in the usual way, 
leading to a semi-discrete problem that is subsequently propagated in time by numerical integration.

When dealing with convection-dominated problems, the spatial discretization can be chosen to evolve continuously in time,
which allows the finite element mesh to continuously track moving structures in the solution such as steep sweeping fronts 
\cite{baines1994moving,CARLSONMILLER1,CARLSONMILLER2}.
These moving finite elements can lead to remarkably improved stability in computing a solution, with respect to the length of permissible time steps
\cite{MILLER1,MILLER2}.

In \cite{BANKMETTI,BANKSANTOS}, tensor-product finite element spaces compatible with a method of lines discretization were
introduced that allowed these moving finite element solutions to be studied in a space-time finite element framework.
As a result, these papers established symmetric error estimates for these finite element solutions
when the numerical time integrator belongs to a particular class of collocation methods.
The first such symmetric error estimate is proven in \cite{DUPONT82} for semi-discrete moving finite elements
by using a mesh-dependent energy semi-norm, $\enorm{\cdot}$. 
To elaborate, a symmetric error estimate states that the finite element solution, $u_h$, satisfies
\begin{equation}        \label{symmetric error estimate}
        \enorm{u - u_h}         \le     C \inf_{\chi\in\mathcal{V}_h} \enorm{u - \chi},
\end{equation}
where $u$ is the true solution to the differential equation and $\mathcal{V}_h$ is the (tensor-product) moving finite element space.

In this paper, we consider the effects of employing a time integrator that does not belong to the previous class of collocation methods.
This is a valuable modification because collocation methods implicitly couple all intermediate stages of each time step,
significantly increasing the computational complexity when using higher order quadrature.
We consider the second-order and diagonally-implicit time integrator TR-BDF2 \cite{A30,A30a},
and using piecewise quadratic tensor-product finite element spaces to discretize the problem.
This time integration scheme is known for its favorable stability properties \cite{STRANG,hosea1996analysis,ying2009composite}.
Moving finite element discretizations often lead to stiff systems of ODEs \cite{MILLER1,MILLER2},
which is why a stable time integration scheme is required.
In Section \ref{sec3}, we prove an error estimate like \eref{symmetric error estimate} with an additional term 
introduced by TR-BDF2.

This work largely builds on the analyses in \cite{BANKMETTI,METTITHESIS}, where parts of the preliminary analysis are given in more detail.
This paper is organized as follows: in Section \ref{sec2}, we describe the model equation, the piecewise quadratic tensor-product finite element space,
and some preliminary results.
In Section \ref{sec3}, a space-time moving finite element method using TR-BDF2 time integration is proposed and an error estimate for the finite element solution is proven.
Section \ref{sec4} describes and reports some numerical experiments that validate the efficacy of these moving finite element methods.

\section{Preliminary Results}\label{sec2}

The model problem used in this error analysis is the linear convection-diffusion equation.
The spatial domain, $\Omega$, is assumed to be a simply connected set in $\Re^d$,
where {$d=1,2,$ or $3$}, with boundary $\partial \Omega$.
The time domain is a finite interval, $(0, T]$, and the space-time domain is given by $\mathcal{F} \equiv \Omega \times (0, T]$.

Let $a$, $b$, $c$, and $f$ be smooth and bounded functions defined on $\mathcal{F}$
such that there exist constants $\bar{a}>0$ and $\bar{c}\ge0$ with $a \ge \bar{a}$ and $c \ge \bar{c}$ on $\Omega$,
and let $g$ be piecewise continuous on $\partial \Omega$.
Let $u_0$ be a given initial condition for the solution on $\mathcal{F}$ and let $n$ denote the outward unit normal vector to the boundary $\partial \Omega$.
The solution to the differential equation, denoted by $u$, is the function that satisfies 
\begin{align}
        u_t     -       \nabla \cdot ( a \nabla u )     +       b \cdot \nabla u        +       c u     &=      f,              &\mathrm{in}\ \mathcal{F},
                                                                                                                                                        \label{strong pde}\\
                                                                                a  \nabla u \cdot n     &=      g,              &\mathrm{on}\ \partial \Omega \times (0, T],
                                                                                                                                                        \label{strong bc}\\
                                                                                                u(x, 0) &=      u_0(x), &\mathrm{for}\ x\ \mathrm{in}\ \Omega.\nonumber
\end{align}

When the convection term, $b$ in \eref{strong pde}, is large relative to the other coefficients,
steep shock layers may develop in the solution that propagate through the spatial domain.
Basic finite element discretizations consequently require short time steps to maintain accuracy of the computed solution at this moving front,
or a moving mesh can be employed, allowing for more flexibility in the length of the time step.
Moving finite elements offset the convection velocity in a similar manner to the method of characteristics.
We assume a time dependent parametrization of the spatial domain, $x(t)$, and define the space-time derivative as
\[
        u_\tau(x(t),t)  \equiv  \frac{d}{dt} u(x(t),t)  =       u_t (x(t),t) + x_t \cdot \nabla u(x(t),t).
\]
We refer to this as the \emph{characteristic derivative} of $u$.

The weak form of the problem is:
find $u$ with $u(t) \in \mathcal{H}^1(\Omega)$ and $u_t(t) \in \mathcal{L}_2(\Omega)$
such that for all $\chi$ in $\mathcal{H}^1(\Omega)$ and $0 < t \le T$,
\begin{align}
                        \bigip{ u_\tau(\cdot,t)}{\chi}  +       \mathcal{A}_\tau\ip{t;  u}{\chi}        &=      \bigip{f(\cdot,t)}{\chi}        +       \bigbip{g(\cdot,t)}{\chi},      
                                                                                                                                                                        \label{variational form diff eqn} \\
\intertext{and when $t=0$}
                                                                                \bigip{u(\cdot,0)}{\chi}        &=      \bigip{u_0}{\chi}.      \nonumber
\end{align}
The inner-products are given by
\begin{align*}
        \ip{f}{\chi}                                    &=      \int_{\Omega} f(x) \chi(x) \ dx,                        \\
        \bip{g}{\chi}                           &=      \int_{\partial \Omega} g(s) \chi(s)     \ ds,
\end{align*}
and define the time-dependent bilinear form
\begin{multline*}
        \mathcal{A}_\tau\ip{t; u}{\chi} \equiv  \int_{\Omega}   a(x,t) \nabla u(x,t) \cdot \nabla \chi(x)       +       (b(x,t)-x_t(t)) \cdot \nabla u(x,t) \ \chi(x)   \\
                                                                +       c(x,t) u(x,t) \chi(x)   \ dx.
\end{multline*}
Notice that parameterizing the spatial variable so that $x_t \approx b$ leads to a formulation where the convection velocity is much less prominent,
as it is ``absorbed'' into the characteristic derivative.

The finite element space we use to discretize the differential equation is a tensor-product of a discontinuous piecewise quadratic finite elements in time
with continuous piecewise quadratic finite elements in space.

Let $0 < t_1 < \ldots < t_m = T$ form a strict ordered partition of the time domain and define $\Delta t_i \equiv t_i - t_{i-1}$.
For $1 \le i \le m$, let $\{x_k(t)\}$ represent the vertices of a triangulation of the domain at time $t$, where $0 \le k \le N_i$,
and we assume that $| x_k(t)-x_j(t) | > \Delta x$ throughout each time partition for some minimum mesh size $\Delta x>0$ and $j \neq k$.
The vertices of the mesh are permitted to move along quadratic trajectories throughout each time partition---that is, 
each node $x_k(t)$ is a quadratic polynomial for $t_{i-1} < t \le t_i$ in time---though
discontinuous reconfigurations of the mesh are permitted at the beginning of each partition.
These discontinuous changes in the mesh provide flexibility for periodically adding and removing degrees of freedom,
as well as keeping the nodes in the mesh from colliding and tangling.

\begin{figure}[h!btp] 
        \begin{picture}(300,100)(-35,0)
                \put(0,0){\thicklines{\line(0,1){70}}}
                \put(0,0){\thicklines{\line(1,0){300}}}
                \put(0,35){\dashbox{4}(300,0){}}
                \put(0,70){\thicklines{\line(1,0){300}}}
                \put(300,0){\thicklines{\line(0,1){70}}}
                \put(0,0){\thicklines{\circle*{5}}}
                \put(300,70){\thicklines{\circle*{5}}}
                \put(300,0){\thicklines{\circle*{5}}}
                \put(0,70){\thicklines{\circle*{5}}}
                \put(250,70){\thicklines{\circle*{5}}}
                \put(200,70){\thicklines{\circle*{5}}}
                \put(130,70){\thicklines{\circle*{5}}}
                \put(60,70){\thicklines{\circle*{5}}}
                \put(40,0){\thicklines{\circle*{5}}}
                \put(100,0){\thicklines{\circle*{5}}}
                \put(170,0){\thicklines{\circle*{5}}}
                \put(260,0){\thicklines{\circle*{5}}}
                \put(47,35){\thicklines{\circle{5}}}
                \put(111.5,35){\thicklines{\circle{5}}}
                \put(178,35){\thicklines{\circle{5}}}
                \put(252,35){\thicklines{\circle{5}}}
                \put(0,35){\thicklines{\circle{5}}}
                \put(300,35){\thicklines{\circle{5}}}
                \put(20,0){\circle{4}}
                \put(68,0){\circle{4}}
                \put(133,0){\circle{4}}
                \put(213,0){\circle{4}}
                \put(279,0){\circle{4}}
                \put(21,35){\circle{4}}
                \put(77,35){\circle{4}}
                \put(143,35){\circle{4}}
                \put(214,35){\circle{4}}
                \put(275,35){\circle{4}}
                \put(28,70){\circle{4}}
                \put(90,70){\circle{4}}
                \put(160,70){\circle{4}}
                \put(222,70){\circle{4}}
                \put(273.5,70){\circle{4}}
                \qbezier(40,0)(44,35)(60,70)
                \qbezier(100,0)(108,35)(130,70)
                \qbezier(170,0)(171,35)(200,70)
                \qbezier(260,0)(249,35)(250,70)
        \end{picture}
\caption{An example space-time mesh partition of a single dimension in space.
The filled circles represent the space-time ``hat'' basis nodes;
hollow circles correspond to basis nodes with basis functions that are the product of a ``bump'' function with a ``hat'' or
``bump'' function.}
\label{fig: space time partition}
\end{figure}
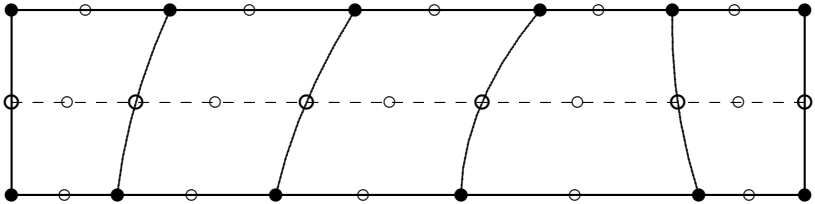


The reference element for this finite element space is the Cartesian product of the unit interval ($d=1$), triangle ($d=2$), or tetrahedron ($d=3$)
for space with the unit interval reference element for the time domain.
Let $e$ be an element in the mesh with vertices given by $x_k(t)$ at time $t$, where $k=0,1,\ldots,d$ and $t_{i-1} \le t \le t_i$.
Then the (isoparametric) mapping from the reference element to $e$ is given by
\begin{align*}
        t       &=      t_{i-1} + \hat{t} \Delta t_i,\\
        x(t)    &=      \mathcal{J}_e(t) \hat{x} + x_{0}(t),
\end{align*}
for $0 \le \hat{t} \le 1$ and $\hat{x}$ in the spatial reference element,
where $x(t)$ is affine in space and quadratic in time, as the node trajectories are restricted to follow quadratic paths throughout the time partition.
Notice that the time variable is space invariant, but the spatial variable does in fact depend on $t$.
The $d \times d$ spatial Jacobian matrix, $\mathcal{J}_e(t)$, determines the shape and size of the element $e(t)$.
Since the determinant of the spatial Jacobian is proportional to the size of the element, 
we require $| \det \mathcal{J}_e(t) | > 0$ for all $t$ to ensure a non-degenerate mesh.

Since we are using a tensor-product space-time discretization, the degrees of freedom of the finite element space are distributed in time slices.
Namely, for fixed time $t$, the degrees of freedom $\{x_k(t)\}$ define a standard finite element space of continuous piecewise-quadratic polynomials on $\Omega$,
which we denote by $\mathcal{V}_h(t)$.
This is an important property that is directly used in formulating the discrete problem.
We emphasize that the finite element functions are piecewise quadratic polynomials along the trajectories of the space nodes $x_k(t)$ 
defined via the isoparametric mapping, rather than in the time direction as with standard tensor-product discretizations.
This ensures that the \emph{characteristic} derivative of a finite element function, $\phi$,
is continuous within each time partition and satisfies $\phi_\tau(t) \in \mathcal{V}_h(t)$.
We denote the tensor-product finite element space on the space-time domain by $\mathcal{V}_h$.
More detailed descriptions of these tensor product finite element spaces can be found in \cite{BANKMETTI,METTITHESIS}.

For the error analysis of our moving finite element method,
it is convenient to define the finite element functions at the mesh discontinuity by an $\mathcal{L}_2$-projection.
For $\phi$ in $\mathcal{V}_h$, we represent the limiting values of $\phi$ near the discontinuities as
\(
        \phi(t_{i^+}) \equiv \lim_{\delta\rightarrow0^+} \phi(t+\delta)
\)
and
\(
        \phi(t_{i^-}) \equiv \lim_{\delta\rightarrow0^+} \phi(t-\delta).
\)
Thus, we require
\[
        \ip{\phi(t_{i^+})}{\chi}        =       \ip{\phi(t_{i^-})}{\chi}
\]
for all $\chi$ in $\mathcal{V}_h(t_{i^+})$,   $i=1,\ldots,m$.
To uniquely determine the finite element functions, we take $\phi(t_i) \equiv \phi(t_{i^-})$ at the discontinuities.

Multi-index notation is used to represent spatial derivatives, 
but time and characteristic derivatives will not follow this convention.
The $\mathcal{H}^k(\Omega)$ semi-norm and norm follow conventional notation and we write
\[
        \snorm{v}_k    =      \left(\sum_{|\alpha| = k} \ip{ D_{\alpha} v }{ D_{\alpha} v } \right)^{1/2}    
\quad\mathrm{and}\quad
        \norm{v}_k      =      \left(\sum_{|\alpha| \le k} \ip{ D_{\alpha} v }{ D_{\alpha} v } \right)^{1/2}.
\]
Following Dupont \cite{DUPONT82}, a mesh-dependent semi-norm is defined that allows us to prove our error estimate,
\[
        \norm{v}_{(-1,\mathcal{V}_h^p(t))}      =       \sup_{\substack{ \chi \in \mathcal{V}_h^p(t) \\ \chi \neq 0}}   \frac{ \left| \ip{ v }{\chi} \right|}{ \norm{\chi}_1}.
\]
We also use the infinity norm, $\norm{v}_{\infty}       =       \max_{x\in\Omega} \snorm{v(x)}$.

We now introduce a space-time shape regularity constraint for the moving finite elements,
that controls the time evolution of the spatial elements 
and prevents degenerate elements.
Fix $e$ to be an element in the time partition with $t_{i-1} \le t \le t_i$.
Then, the Jacobian matrix at time $t$ can be represented as
\begin{equation}        \label{jacobian evolution}
        \mathcal{J}_e(t)        =       \big( \mathcal{R}_{e}(t) + \Delta t_i \mathcal{H}_e(t) \big) \mathcal{J}_e( t_{i-1^+} ),
\end{equation}
for some orthogonal \emph{rotation} matrix, $\mathcal{R}_{e}(t)$, and \emph{evolution} matrix, $\mathcal{H}_e(t)$.
The matrix $\mathcal{R}_{e}+\Delta t_i \mathcal{H}_e$ is constrained to have quadratic polynomial entries throughout the time partition.
The matrix $\mathcal{R}_{e}(t)$ describes the element rotation in time,
 and the evolution matrix describes the deformation of the shape of the element.
Since the trajectories of the spatial nodes are restricted to quadratic polynomial paths,
elements cannot rotate perfectly in time and more of a twisting action is observed;
the evolution matrix necessarily reflects these deformations.
If an element is merely translated in time, without rotation or changing shape,
then the Jacobian matrix, $\mathcal{J}_e(t)$, remains unchanged.

Let $\rho(\cdot)$ represent the spectral radius for $d \times d$ matrices.
It is assumed that the evolution matrix $\mathcal{H}_e$ has a uniformly bounded spectral radius throughout the time step;
namely, there exists some positive constant $\mu$ that does not depend on $e$ or $t$ such that
\begin{equation}        \label{space time regularity}
        \rho \big( \mathcal{H}_e(t) \big)       \le \mu.
\end{equation}
This bounds the relative change in shape and size of the element over time.
Assuming a non-degenerate finite element space and the space-time shape regularity bound \eref{space time regularity}, it follows that
\begin{equation}        \label{spectral radius bound}
        \rho \big( \mathcal{J}_e(t) \mathcal{J}_e^{-1}(t_{i-1^+}) \big) =       \rho \big( \mathcal{R}_{e}(t) + \Delta t_i \mathcal{H}_e(t) \big)
                                                                                        \le     1 + \mu \Delta t_i
\end{equation}
and, for $\tilde{c}_{\mu,d} = [ (1+\mu \Delta t_i)^d - 1] / \Delta t_i = \mathcal{O}(1)$ and $\Delta t_i \le 1/2\tilde{c}_{\mu,d}$,
\begin{multline}        \label{determinant bound}
        (1-\tilde{c}_{\mu,d}\Delta t_i) \le     (1- \mu \Delta t_i)^d   \le     \frac{\mathcal{D}_e(t)}{\mathcal{D}_e(t_{i-1^+})}               \\
                        =       \det \big(\mathcal{J}_e(t) \mathcal{J}_e^{-1}(t_{i-1^+}) \big)  \le     (1 + \mu \Delta t_i)^d  \le     (1+ \tilde{c}_{\mu,d}\Delta t_i),
\end{multline}
since $\mathcal{R}_{e}$ is an orthogonal matrix.

Let $\phi$ be a function in the finite element space $\mathcal{V}_h(t)$ for some $t$ in the time partition $(t_{i-1},t_i]$.
We \emph{shift} $\phi$ onto the mesh of $\mathcal{V}_h^p(t_{i-1^+})$, at the beginning of the time partition by replacing
the basis functions of $\mathcal{V}_h(t)$ with with their corresponding basis functions in $\mathcal{V}_h(t_{i-1^+})$,
while preserving the basis coefficients.
Formally, this operation can be defined by an element-wise composition of the inverse of the affine spatial isoparametric maps for the elements in the mesh at time $t$,
which is well-defined for non-degenerate meshes,
with the affine spatial isoparametric maps for the elements at the beginning of the time step $t_{i-1^+}$.
The following lemma, proven in \cite{METTITHESIS}, establishes the relationship between the space-time shape regularity constraint \eref{space time regularity}
and the continuity of this shift operation.

\begin{lemma}[Shift Lemma]      \label{shift lemma}
Let $\phi,\chi \in \mathcal{V}_h(t)$ and $\tilde{\phi},\tilde{\chi} \in \mathcal{V}_h(t_{i-1^+})$ represent a pair of finite element functions and their shifts, respectively,
on a non-degenerate time partition of the mesh that satisfies \eref{space time regularity} on each element.
If $\Delta t_i \le 1/2\tilde{c}_{\mu,d}$, as defined in \eref{determinant bound}, then there exists a positive constant $C_{\mu,d}$ such that
\begin{align}
        \Big| \bigip{\phi}{\chi} - \bigip{\tilde{\phi}}{\tilde{\chi}} \Big|     &\le    C_{\mu,d} \Delta t_i \frac{\bignorm{\tilde{\phi}}_0^2   +       \bignorm{\tilde{\chi}}_0^2}{2}, \label{IP1}\\
        \Big|   \bignorm{\phi}_0^2 - \bignorm{\tilde{\phi}}_0^2 \Big|   &\le    C_{\mu,d} \Delta t_i \bignorm{\tilde{\phi}}_0^2,                                                                          \label{L2 bound}\\
        \Big|   \bignorm{\phi}_1^2 - \bignorm{\tilde{\phi}}_1^2 \Big|   &\le    C_{\mu,d} \Delta t_i \bignorm{\tilde{\phi}}_1^2.                                                                          \label{H1 bound}
\end{align}
\end{lemma}

We now present a local Gr\"onwall lemma that will be used to bound the maximum error of the finite element solution in the $\mathcal{L}_2$-norm over each time partition.
The original proof for this lemma is given in \cite{METTITHESIS}.

\begin{lemma}[Local Gr\"onwall Inequality]      \label{local gronwall}
Suppose there are two distinct times on each time partition, $t_{i-1} < t_{i,1} < t_{i,2} \le t_i$,
when the mesh satisfies the regularity constraint \eref{space time regularity} and that there exists a positive constant $\kappa$ such that
\begin{equation}        \label{mesh alignment local lemma}
        \norm{b - x_t}_\infty   \le     \kappa.
\end{equation}
If $\Delta t_i \le1/2\tilde{c}_{\mu,d}$, as defined in \eref{determinant bound},
and functions $\phi$ in $\mathcal{V}_h$ and $\eta$ with $\eta$ in $\mathcal{H}^1(\Omega)$ and $\eta_{\tau}$ in $\mathcal{L}_2(\Omega)$ satisfy
\begin{equation}        \label{collocation constraint}
        \bigip{\phi_{\tau}(t_{i,j})}{\chi}      +       \mathcal{A}_\tau\bigip{\phi(t_{i,j})}{\chi}      
                        =       \bigip{\eta_{\tau}(t_{i,j})}{\chi}      +       \mathcal{A}_\tau\bigip{\eta(t_{i,j})}{\chi}      
\end{equation}
for all $\chi$ in $\mathcal{V}_h(t_{i,j})$ at time each collocation node $j=1,2$,
then, there exists a constant such that
\begin{multline*}
        \max_{1 \le j \le 2} \bignorm{ \phi(t_{i,j}) }_0^2    \\  \le
                C \left\{       \bignorm{ \phi(t_{i-1^+})}_0^2 
                + \sum_{j=1}^2 \Delta t_i\left( \norm{\eta_{\tau}(t_{i,j})}_{(-1,\mathcal{V}_h(t_{i,j}))}^2 + \norm{\eta(t_{i,j})}_1^2 + \norm{\phi(t_{i,j})}_1^2 \right) \right\},
\end{multline*}
where $C$ depends on $\kappa, \mu, d, p$, and the differential equation.
\end{lemma}

Another discrete Gr\"onwall lemma is used to aggregate the spatial error bounds from each time partition over the entire time domain.
This result was proven in \cite{SANTOSTHESIS}.

\begin{lemma}[Discrete Gr\"onwall Inequality]   \label{discrete gronwall}
Let $\Delta t_i > 0$ and $\alpha_i, \gamma_i, \theta_i, q_i \ge 0$, for $1 \le i \le m$, with $\theta_i\Delta t_i \le \frac{1}{2}$ and $\theta = \max_i \theta_i$.
Then, if
\[
        \frac{q_i - q_{i-1}}{\Delta t_i} + \gamma_i     \le     \alpha_i + \theta_i ( q_i + q_{i-1} ),
\]
there exists a positive constant $C_{\theta}$ such that
\[
        \max_{1 \le i \le m} q_i + \sum_{i=1}^m \gamma_i \Delta t_i             \le C_{\theta} \left\{ q_0 + \sum_{i=1}^m \alpha_i \Delta t_i \right\}.
\]
\end{lemma}

\section{A space-time moving finite element method with TR-BDF2} \label{sec3}

To achieve second-order accuracy for the finite element solution, we employ the two-stage, diagonally-implicit time integration scheme TR-BDF2.
Since this scheme is diagonally implicit, significant time savings are realized when it is applied to large systems of ODEs,
including systems arising from a method of lines discretization of a parabolic equation.
TR-BDF2 was proposed by Bank, et al.~in \cite{A30,A30a} and has been analyzed in several other papers
for its efficiency and stability \cite{STRANG,hosea1996analysis,ying2009composite}.
TR-BDF2 actually refers to a family of time stepping methods that is parametrized by the location of the intermediate basis node.
For this method, we define the collocation nodes on the reference element to be 
\[
        \hat{t}_0 = 0,  \quad   \hat{t}_1 = {\varepsilon}/{2},  \quad\mathrm{and}\quad  \hat{t}_2 = 1,
\]
where $0 < \varepsilon <1$ is a free parameter that determines the exact time stepping scheme within the TR-BDF2 family.
We choose the reference element's basis nodes to be
\[
        \hat{\zeta}_{0} = 0,    \quad   \hat{\zeta}_{1} = {\varepsilon},        \quad\mathrm{and}\quad  \hat{\zeta}_{2} = 1
\]
so that $\hat{t}_1$ is the midpoint of the first two basis nodes.
The Runge-Kutta coefficients correspond to integrating the computed solution a step of length $\varepsilon\Delta t_i$ by the trapezoid rule,
then completing the time step by a second-order backward Euler difference:
let $u(t)$ denote an approximate solution of some system of ODEs, then the solution is computed using the approximations
\begin{align*}
        u_{TR}(t_{i,1}) &=      \frac12 u( \zeta_{i,1} ) + \frac12 u( \zeta_{i,0} ),    \\
        u(t_{i,2})              &=      u( \zeta_{i,2} ),
\end{align*}
at the collocation nodes, where the isoparametric map is use to distribute the basis and collocation nodes throughout the time domain
\[
        \zeta_{i,j}     = t_{i-1} + \hat{\zeta}_j \Delta t_i 
\quad\mathrm{and}\quad
        t_{i,k}         = t_{i-1} + \hat{t}_k \Delta t_i
\]
for $j=0,1,2$ and $k=1,2$.
The coefficients for the time derivative are determined by the interpolating quadratic Lagrange polynomials associated
with the degrees of freedom in the time discretization evaluated at the collocation nodes:
\begin{align*}
        \bar{\delta}_t u (t_{i,1})      &=      \frac{1}{\varepsilon} u( \zeta_{i,1} ) -\frac{1}{\varepsilon} u (\zeta_{i,0}), \\
        \bar{\delta}_t u (t_{i,2})      &=      \frac{2-\varepsilon}{1-\varepsilon}u( \zeta_{i,2} )
                                                        -        \frac{1}{\varepsilon(1-\varepsilon)}u( \zeta_{i,1} )
                                                        +       \frac{1-\varepsilon}{\varepsilon}u(\zeta_{i,0}).
\end{align*}

The optimal choice for the parameter is known to be $\varepsilon = 2 - \sqrt2$, as it minimizes the local truncation error 
and we refer to the choice $\varepsilon = 2- \sqrt2$ as the Richardson basis node \cite{choudhury1992waveform,ying2009composite}.
The TR-BDF2 scheme is A-stable and L-stable \cite{A30,A30a,hosea1996analysis}.
These stability properties are important as method of lines discretizations of \eref{variational form diff eqn} typically lead to stiff systems of ODEs
\cite{MILLER1,MILLER2}.

In this section, we use the TR-BDF2 method to integrate the semi-discrete system of ODEs given by the method of lines applied to a 
moving finite element discretization, as described in Section \ref{sec2} of the weak form of the 
differential equation \eref{variational form diff eqn}.
Fix $i$ to index some time partition and let $u_h \in \mathcal{V}_h$;
let $\tilde{u}_h(t)$ represent the function shifted onto the mesh at the mid-step collocation node $t=t_{i,1}$
and ${\hat{u}}_h(t)$ represent the function shifted onto the mesh at the end-step collocation node $t=t_{i,2}$.
The discrete problem can be characterized as finding the finite element function ${u}_h$ that satisfies
\begin{equation}        \label{midstep constraint}
        \Bigip{\frac{\tilde{u}_h(\zeta_{i,1})-\tilde{u}_h(\zeta_{i,0})}{\varepsilon \Delta t_i}}{\tilde{\chi}}
                +       \mathcal{A}_{\tau}\Bigip{\frac{\tilde{u}_h(\zeta_{i,1})+\tilde{u}_h(\zeta_{i,0})}{2}}{\tilde{\chi}}
                                                =       \bigip{f(t_{i,1})}{\tilde{\chi}}        +       \bigbip{g(t_{i,1})}{\tilde{\chi}}
\end{equation}
for all $\tilde{\chi}$ in $\mathcal{V}_h(t_{i,1})$ at each mid-step collocation node $t=t_{i,1}$,
and for all ${\hat{\chi}}$ in $\mathcal{V}_h(t_{i,2})$,
\begin{multline}        \label{endstep constraint}
        \Bigip{\frac{   \varepsilon(2-\varepsilon) {\hat{u}}_h(\zeta_{i,2}) -   {\hat{u}}_h(\zeta_{i,1}) + (1-\varepsilon)^2 {\hat{u}}_h(\zeta_{i,0})}
                                {\varepsilon(1-\varepsilon) \Delta t_i}}{{\hat{\chi}}}  +       \mathcal{A}_{\tau}\ip{{\hat{u}}_h(\zeta_{i,2})}{{\hat{\chi}}}
                                        \\      =       \bigip{f(t_{i,2})}{{\hat{\chi}}}        +       \bigbip{g(t_{i,2})}{{\hat{\chi}}},
\end{multline}
$i=1,\ldots,m$.
In contrast to the weak formulation \eref{variational form diff eqn}, the discrete problem weakly enforces the differential equation only at the collocation nodes.
The constraint
\begin{equation}        \label{jump condition}
        \ip{{u}_h(t_{i^+})}{\chi}       =       \ip{{u}_h(t_{i^-})}{\chi}
\end{equation}
must also hold for all $\chi$ in $\mathcal{V}_h(t_{i^+})$, $i=1,\ldots,m$, to ensure that ${u}_h \in \mathcal{V}_h$.

We now prove an error estimate for \eref{midstep constraint}--\eref{jump condition}.
Our proof follows that of Theorem 4.3 in \cite{BANKMETTI} with some additional arguments
that bound the error introduced by the trapezoid approximation at the mid-step of each time partition.
Due to the departure of this method from a strict finite element framework, the symmetry of the error bound is broken and an additional term
proportional to the error of the trapezoid approximation is introduced.
We must also assume the bound
$\norm{ a_{{\tau}} }_{\infty}    \le     \alpha$, 
for some $\alpha > 0$  on the characteristic derivative of the diffusion coefficient.

One final aspect in which the error bound changes is that the TR-BDF2 increases sensitivity to the discontinuous changes in the mesh at the beginning of the time steps,
compared to the collocation methods used in \cite{BANKMETTI}.
Recently, Bank and Yserentant \cite{BANKYSERENTANT} proved the $\mathcal{H}^1$-stability of $\mathcal{L}_2$-projections
onto finite element spaces with potentially nonuniform meshes.
Using this result, we assume $\snorm{\chi(t_{i^+})}_1 \le C_{\mathcal{H}} \snorm{\chi(t_{i^-})}_1$, for $\chi$ in $\mathcal{V}_h$.
As can be seen from \cite{BANKYSERENTANT}, the bounding constant, $C_{\mathcal{H}}$, 
is smaller when the mesh reconfiguration is more subtle at the mesh discontinuities.
This intuitively makes sense, since $\chi(t_{i^+}) \approx \chi(t_{i^-})$ in such cases.
For a given differential equation, since $0 < \bar{a} \le a \le \norm{a}_{\infty}$, we have an equivalence of norms
\[
        \frac{1}{c_{\mathcal{A}}} \norm{a^{1/2}(t_i) \nabla \phi}_0 \le \snorm{\phi}_1 \le c_{\mathcal{A}} \norm{a^{1/2}(t_i) \nabla \phi}_0
\]
for some positive $c_{\mathcal{A}}$.
We define the stability constant for the diffusion-weighted semi-norm by
\begin{equation}
        \snorm{a^{1/2}(t_i)\nabla\phi(t_{i^+})}_0^2     \le     C_{\mathcal{A},\mathcal{H}} \snorm{a^{1/2}(t_i)\nabla\phi(t_{i^-})}_0^2.
\end{equation}

The norm in which the error is bounded employs the trapezoid approximation at the mid-step collocation node of each time partition.
Let $u$ be some function defined on $\mathcal{F}$ and $u_{TR}(t_{i,1}) \equiv (\tilde{u}(\zeta_{i,1}) + \tilde{u}(\zeta_{i,0}))/{2} \in \mathcal{V}_h(t)$,
where $\tilde{u}(t)$ represents $u$ evaluated at time $t$ following the characteristic from time $t_{i,1}$
(meaning $\tilde{u}(t)$ represents the function $u$ shifted onto the mesh at time $t_{i,1}$.)
The semi-norm in which we bound the error of the finite element solution with TR-BDF2 time integration is
\begin{multline*}
        \enorm{u}^2     \equiv  \max_{\substack{1 \le i \le m \\ 1 \le j \le 2}} \norm{u(t_{i,j})}_0^2
                                                        +       \sum_{i=1}^m \Delta t_i \big(   \norm{ \bar{\delta}_{{\tau}} u(t_{i,1}) }_{(-1,\mathcal{V}_h(t_{i,1}))}^2 + \norm{ u_{TR}(t_{i,1}) }_1^2 
                                                        \\+                                                     \norm{ u(\zeta_{i,1}) }_1^2
                                                        +                                                       \norm{ \bar{\delta}_{{\tau}} u(t_{i^-}) }_{(-1,\mathcal{V}_h(t_{i^-}))}^2 + \norm{ u(t_{i^-}) }_1^2     \big).
\end{multline*}

\begin{theorem} \label{trbdf a priori theorem}
Suppose that $\mathcal{V}_h$ is a finite element space with a non-degenerate mesh and let $\varepsilon = 2/3$.
Furthermore, assume that there exist positive constants $\alpha$, $\mu$, and $\kappa$ such that at each collocation node
\begin{align}
        \norm{a_{\tau}}_{\infty}                        &\le    \alpha, \label{diffusion bound}                 \\
        \rho \big(\mathcal{H}_e\big)         &\le    \mu,            \label{spacetime regularity trbdf}      \\
        \norm{b - x_t}_\infty                           &\le    \kappa, \label{characteristic regularity trbdf}
\end{align}
and that the mesh discontinuities are controlled and the spatial meshes and length of the time steps are graded so that
\begin{align}
        C_{\mathcal{A},\mathcal{H}}             &\le    {7}/{2}                         \label{stability bound} \\
\intertext{and}
        \Delta t_i                                              &\le    2 \Delta t_{i-1}.               \label{graded time step}
\end{align}
Then, if $\Delta t = \max_{1\le i \le m} \Delta t_i$ is sufficiently small, there exists a positive constant $C$ such
that the finite element solution satisfies
\begin{equation}        \label{a priori trbdf}
        \enorm{u - \bar{u}}^2   \le     C \Big\{        \inf_{\chi \in \mathcal{V}_h} \enorm{u - \chi}^2
                                                        +       \int_0^T \norm{ \Delta t_i^2 u  _{\tau\tau}(t) }_1^2\ dt        \Big\},
\end{equation}
where $C$ depends on $\mu, \kappa, d, p$, and the differential equation.
\end{theorem}

Note that this proof is restricted to the case where we use the collocation nodes determined by Gauss-Radau quadrature, $\varepsilon = 2/3$,
as Gauss-Radau quadrature has a positive truncation error, which helps bound the aggregation of the local truncation errors,
and fixes a collocation node at the end of the time steps, which is required for the TR-BDF2 scheme.
Notice that the intermediate time basis node for Gauss-Radau is $2/3$, which is close to the optimal value $\varepsilon = 2-\sqrt{2}\approx0.5858$ for TR-BDF2.
Also, the assumptions \eref{stability bound} and \eref{graded time step} simplify the proof,
although they are stricter than necessary.

\begin{proof}
For this proof, we use the discrete Galerkin orthogonalities
\[
        \ip{u_\tau-{u}_{h,\tau}}{\chi}  +       \mathcal{A}_\tau\ip{u-{u}_{h,TR}}{\chi} =       0,
\]
at $t=t_{i,1}$ for $\chi$ in $\mathcal{V}_h(t_{i,1})$ and
\[
        \ip{u_\tau-{u}_{h,\tau}}{\hat{\chi}}  +       \mathcal{A}_\tau\ip{u-{u}_{h}}{\hat{\chi}} =       0,
\]
at $t=t_{i^-}$ for $\hat{\chi}$ in $\mathcal{V}_h(t_{i^-})$, for $i=1,\ldots,m$.
Following Dupont \cite{DUPONT82}, let $\psi$ in $\mathcal{V}_h$ be an arbitrary function and define 
$\phi \equiv {u}_h - \psi$ in $\mathcal{V}_h$ and $\eta \equiv u - \psi$.
Then, we have
\begin{multline}
        \ip{\phi_\tau}{\chi}    +       \mathcal{A}_\tau\ip{\phi_{i,TR}(t_{i,1})}{\chi} 
                        =       \ip{\eta_\tau}{\chi}    +       \mathcal{A}_\tau\ip{\eta_{i,TR}(t_{i,1})}{\chi}
                                \\        +       \mathcal{A}_\tau\ip{u(t_{i,1}) - u_{i,TR}(t_{i,1})}{\chi}       \label{tr step}
\end{multline}
and
\begin{equation}
        \ip{\phi_\tau}{\hat{\chi}}      +       \mathcal{A}_\tau\ip{\phi(t_{i,2})}{\hat{\chi}(t_{i,2})}
                        =       \ip{\eta_\tau}{\hat{\chi}}      +       \mathcal{A}_\tau\ip{\eta(t_{i,2})}{\hat{\chi}}.          \label{trbdf end step}  
\end{equation}

We begin with bounding $\phi$ at the end step \eref{trbdf end step}, which follows the proof of Theorem 4.3 in \cite{BANKMETTI} exactly.
We choose $\hat{\chi} = \phi(t_{i^-})$ to get
\begin{align}
        \ip{\phi_\tau}{\phi}                                            &=      \frac12 \partial_\tau \norm{\phi(t_{i^-})}_0^2, \label{trbdf first term bound}  \\
        \mathcal{A}_\tau\ip{\phi(t_{i,2})}{\phi(t_{i,2})}       &\ge    (1-\epsilon) \norm{a^{1/2}(t_i) \nabla\phi(t_{i^-})}_0^2        -       C_{\kappa} \norm{\phi(t_{i^-})}_0^2,    \\
        \ip{\eta_\tau}{\phi}                                            &\le    C \Big(\norm{\eta_{\tau}}_{(-1,\mathcal{V}_h(t_{i^-}))}^2 + \norm{\phi(t_{i^-})}_0^2 \Big)
                                                                                        +       \epsilon \norm{a^{1/2}(t_i) \nabla\phi(t_{i^-})}_0^2,   \label{rhs trbdf i}     \\
        \mathcal{A}_\tau\ip{\eta(t_{i,2})}{\phi(t_{i,2})}       &\le    C_{\kappa} \Big( \norm{\eta}_1^2 + \norm{\phi(t_{i^-})}_0^2 \Big)
                                                                                        +       \epsilon \norm{a^{1/2}(t_i) \nabla\phi(t_{i^-})}_0^2,   \label{rhs trbdf ii}
\end{align}
where $\epsilon$ is an arbitrarily small positive constant.
Combining \eref{trbdf first term bound}--\eref{rhs trbdf ii} gives the bound
\begin{multline}        \label{trbdf end step bound}
        \frac12\partial_\tau\norm{\phi(t_{i^-})}_0^2            +       (1-\epsilon)\norm{a^{1/2}(t_i) \nabla\phi(t_{i^-})}_0^2
                                \\\le   C \big\{        \norm{\eta_\tau(t_{i^-})}_{(-1,\mathcal{V}_h(t_{i^-}))}^2       +       \norm{\eta(t_{i^-})}_1^2        +       \norm{\phi(t_{i^-})}_0^2        \big\}.
\end{multline}

The bound at the mid-step collocation node follows from a new choice of test function;
define $\chi$ in \eref{tr step} to be 
\[
        \chi = \phi_{i,TR} + \frac{\Delta t_i}{24} \phi_\tau(t_{i,1})   = \phi(t_{i,1}) + \frac{\Delta t_i}{24} \phi_{\tau}(t_{i,1}) + \frac{\Delta t_i^2}{18}\phi_{i,\tau\tau},
\]
where $\frac{\Delta t_i^2}{18}\phi_{i,\tau\tau}$ is the error of the trapezoid approximation.
The term $\frac{1}{24}\Delta t_i \phi_{\tau}(t_{i,1})$ is not an intuitive addition to the test function;
however, it will lend itself to perfectly offset the additional errors that arise front the trapezoid approximation of ${u}_h$.
Using this test function, we have
\begin{align}
        \lefteqn{\bigip{\phi_\tau(t_{i,1})}{\phi(t_{i,1}) + \frac{\Delta t_i}{24} \phi_{\tau}(t_{i,1}) + \frac{\Delta t_i^2}{18}\phi_{i,\tau\tau}}}             \nonumber\\
                                        &\ge    \frac12 \partial_\tau\norm{\phi(t_{i,1})}_0^2   +       \frac{\Delta t_i}{24}\norm{\phi_{\tau}(t_{i,1})}_0^2    
                                                                -       \frac{\Delta t_i^2}{18}\bignorm{\phi_\tau(t_{i,1})}_0\bignorm{\phi_{i,\tau\tau}}_0              \nonumber\\
                                        &\ge    \frac12 \partial_\tau\norm{\phi(t_{i,1})}_0^2   +       \frac{\Delta t_i}{24}\norm{\phi_{\tau}(t_{i,1})}_0^2
                                                                -       \bigg(\frac{\Delta t_i}{24}\norm{\phi_{\tau}(t_{i,1})}_0^2
                                                                        + \frac{\Delta t_i^{-1}}{54}\norm{\Delta t_i^2\phi_{i,\tau\tau}}_0^2 \bigg)                             \nonumber\\
                                        &\ge    \frac12 \partial_\tau\norm{\phi(t_{i,1})}_0^2   - \frac{\Delta t_i^{-1}}{54}\norm{\Delta t_i^2\phi_{i,\tau\tau}}_0^2.   \label{tr first term bound}
\end{align}
For the bilinear form, we bound this in two parts:
\begin{equation}
        \mathcal{A}_\tau\ip{\phi_{i,TR}}{\phi_{i,TR}} \ge       (1-\delta') \norm{a^{1/2}(t_{i,1}) \nabla\phi_{i,TR}}_0^2
                                                                                -       C_{\kappa} \norm{\phi_{i,TR}}_0^2
\end{equation}  
and
\begin{align}
        \lefteqn{\mathcal{A}_\tau\ip{\phi_{i,TR}}{\Delta t_i\phi_{\tau}(t_{i,1})/24}}                                                                   \nonumber\\
                        &=              \frac{1}{24} \mathcal{A}_\tau\biggip{\frac{\tilde{\phi}(\zeta_{i,1})
                                        + \tilde{\phi}(\zeta_{i-1^+})}{2}}{3\frac{\tilde{\phi}(\zeta_{i,1}) - \tilde{\phi}(\zeta_{i-1^+})}{2}}  \nonumber\\                     
                        &\ge            \frac{1}{32} \Big\{ (1-\tilde{\delta}')\norm{a^{1/2}(t_{i,1}) \nabla\tilde{\phi}(\zeta_{i,1})}_0^2
                                        -       (1+\tilde{\delta})\norm{a^{1/2}(t_{i,1}) \nabla\tilde{\phi}(t_{i-1^+})}_0^2 \Big\}
                                        -       C_{\kappa} \norm{\Delta t_i\phi_{\tau}}_0^2                                                                     \nonumber\\
                        &\ge            \frac{1-\tilde{\delta}'}{32}\norm{a^{1/2}(t_{i,1}) \nabla\tilde{\phi}(\zeta_{i,1})}_0^2
                                        -\frac{1+\tilde{\delta}}{32}\norm{a^{1/2}(t_{i,1}) \nabla\tilde{\phi}(t_{i-1^+})}_0^2 - C_{\kappa} \norm{\Delta t_i\phi_{\tau}}_0^2,\nonumber\\
\end{align}
where $\delta',\tilde{\delta},\tilde{\delta}' > 0$ can be made arbitrarily small at the expense of growing $C_\kappa$.
Bounds on the right side follow the usual arguments.
For arbitrarily small $\epsilon>0$, we use the boundedness of the diffusion coefficient to show

\begin{multline}
        \ip{\eta_\tau}{\phi_{i,TR} + \Delta t_i\phi_{\tau}/24}  \le     C \Big( \norm{\eta_{\tau}(t_{i,1})}_{(-1,\mathcal{V}_h(t_{i,1})}^2
                                +       \bignorm{\phi_{i,TR}}_0^2 + \bignorm{\Delta t_i\phi_{\tau}(t_{i,1})}_0^2        \Big)
                \\              +       \epsilon \big\{ \norm{a^{1/2}(t_{i,1}) \nabla\phi_{i,TR}}_0^2   + \norm{a^{1/2}(t_{i,1}) \nabla\tilde{\phi}(\zeta_{i,1})}_0^2
                                                + \norm{a^{1/2}(t_{i}) \nabla\tilde{\phi}(t_{i^-})}_0^2 \big\}, \label{rhs tr i}
\end{multline}
\begin{multline}
        \mathcal{A}_\tau\ip{\eta_{i,TR}(t_{i,1})}{\phi_{i,TR}(t_{i,1}) + \Delta t_i\phi_{\tau}/24}      \le     C_{\kappa} \Big( \norm{\eta_{i,TR}}_1^2
                        +       \bignorm{\phi_{i,TR}}_0^2 + \bignorm{\Delta t_i\phi_{\tau}(t_{i,1})}_0^2        \Big)   \\
                        +       \epsilon \big\{ \norm{a^{1/2}(t_{i,1}) \nabla\phi_{i,TR}}_0^2   + \norm{a^{1/2}(t_{i,1}) \nabla\tilde{\phi}(\zeta_{i,1})}_0^2
                                                + \norm{a^{1/2}(t_{i}) \nabla\tilde{\phi}(t_{i^-})}_0^2 \big\}, \label{rhs tr ii}
\end{multline}
and
\begin{multline}
        \mathcal{A}_\tau\ip{u(t_{i,1}) - u_{i,TR}}{\phi_{i,TR} + \Delta t_i\phi_{\tau}/24}              \\
                        \le     C_{\kappa} \Big( \norm{u(t_{i,1}) - u_{i,TR}}_1^2       
                                +       \bignorm{\phi_{i,TR}}_0^2 + \bignorm{\Delta t_i\phi_{\tau}(t_{i,1})}_0^2        \Big)   \\
                                +       \epsilon \big\{ \norm{a^{1/2}(t_{i,1}) \nabla\phi_{i,TR}}_0^2   + \norm{a^{1/2}(t_{i,1}) \nabla\tilde{\phi}(\zeta_{i,1})}_0^2
                                + \norm{a^{1/2}(t_{i}) \nabla\tilde{\phi}(t_{i^-})}_0^2 \big\}. \label{rhs tr iii}
\end{multline}

Hence, at the mid-step collocation node corresponding to the trapezoid step of the TR-BDF2 rule, the following bound follows from
\eref{tr first term bound}--\eref{rhs tr iii}: for small $\epsilon, \tilde{\epsilon}, \tilde{\epsilon}'>0$,
\begin{multline}        \label{mid step tr bound}
        \frac12 \partial_\tau\norm{\phi(t_{i,1})}_0^2   -       \frac{\Delta t_i^{-1}}{54}\norm{\Delta t_i^2\phi_{i,\tau\tau}}_0^2              
                                +       (1-\epsilon') \norm{a^{1/2}(t_{i,1}) \nabla\phi_{i,TR}}_0^2                             \\      
                                +       \frac{1-\tilde{\epsilon}'}{32}\norm{a^{1/2}(t_{i,1}) \nabla\tilde{\phi}(\zeta_{i,1})}_0^2
                                -       \frac{1+\tilde{\epsilon}}{32}\norm{a^{1/2}(t_{i,1}) \nabla\tilde{\phi}(t_{i-1^+})}_0^2
                                -       \epsilon \norm{a^{1/2}(t_{i}) \nabla\tilde{\phi}(t_{i^-})}_0^2  \\
                                \le     C       \Big\{  \norm{\eta_\tau(t_{i,1})}_{(-1,\mathcal{V}_h(t_{i,1})}^2        +       \norm{\eta_{i,TR}}_1^2
                                                        +       \norm{u(t_{i,1}) - u_{i,TR}}_1^2        \\
                                                        +       \bignorm{\phi_{i,TR}}_0^2 + \bignorm{\Delta t_i\phi_{\tau}(t_{i,1})}_0^2        \Big\}.
\end{multline}

Thus, we have the spatial bounds at each collocation node.
We now turn to bound the error introduced by the time stepping scheme.
Since we are using Gauss-Radau quadrature, we have
\[
        \frac34f(1/3) + \frac14f(1)     =       \int_0^1 f(t)\, dt      +       \frac{1}{6^3} f'''(\zeta),
\]
for some $\zeta$ in $[0,1]$ and any bounded function $f$ on $[0,1]$.
For the characteristic derivative terms in \eref{trbdf end step bound} and \eref{mid step tr bound},
applying lemma \ref{shift lemma} and this quadrature rule over the time partition gives
\begin{align*}
        \lefteqn{\frac{3\Delta t_i}{4} \partial_\tau\norm{\phi(t_{i,1})}_0^2    +       \frac{\Delta t_i}{4} \partial_\tau\norm{\phi(t_{i,2})}_0^2}                             \\
                &\ge    \frac{3\Delta t_i}{4} \frac{d}{dt}\norm{\bar{\phi}(t_{i,1})}_0^2 + \frac{\Delta t_i}{4} \frac{d}{dt}\norm{\bar{\phi}(t_{i,2})}_0^2          \\
                &\qquad -       C_{\mu,d}\Delta t_i \Big\{ \norm{\Delta t_i\bar{\phi}_t(t_{i,1})}_0 \norm{\bar{\phi}(t_{i,1})}_0 
                                +       \norm{\Delta t_i\bar{\phi}_t(t_{i,2})}_0 \norm{\bar{\phi}(t_{i,2})}_0 \Big\}                                                                        \\
                &=      \norm{\bar{\phi}(t_{i^-})}_0^2        -       \norm{\bar{\phi}(t_{i-1^+})}_0^2
                        +       \frac{\Delta t_i^4}{6^3} \left(\frac{d}{dt}\right)^4 \norm{\bar{\phi}}_0^2
                        -       C \Delta t_i \max_{0 \le j \le 2} \norm{\bar{\phi}(t_{i,j})}_0^2                                                                                              \\
                &\ge    \norm{\phi(t_{i^-})}_0^2        -       \norm{\phi(t_{i-1^+})}_0^2              +       \frac{\Delta t_i^4}{6^3} \partial_\tau^4 \norm{\phi}_0^2
                                                                                                                -       \hat{C}\Delta t_i \max_{0 \le j \le 2} \norm{\phi(t_{i,j})}_0^2,
\end{align*}
where $\bar{\phi}$ denotes the function $\phi$ shifted onto the mesh at the beginning of the time partition, $\mathcal{V}_h(t_{i-1^+})$.
Since $\partial_\tau^3\phi \equiv 0$, we have $\Delta t_i^4 \partial_\tau^4\norm{\phi}_0^2 = 6 \norm{\Delta t_i^2\phi_{i,\tau\tau}}_0^2$.
Hence, applying the quadrature rule to the characteristic derivative gives
\begin{multline}        \label{time term quadrature}
        \frac{3\Delta t_i}{4} \partial_\tau\norm{\phi(t_{i,1})}_0^2     +       \frac{\Delta t_i}{4} \partial_\tau\norm{\phi(t_{i,2})}_0^2
                \\\ge   \norm{\phi(t_{i^-})}_0^2        -       \norm{\phi(t_{i-1^-})}_0^2              +       \frac{1}{36}\norm{\Delta t_i^2\phi_{i,\tau\tau}}_0^2
                                                                                                                -       \hat{C}\Delta t_i \max_{0 \le j \le 2} \norm{\phi(t_{i,j})}_0^2,
\end{multline}
where we used $\norm{\phi(t_{i-1^-})}_0^2 \ge \norm{\phi(t_{i-1^+})}_0^2$ from \eref{jump condition}.
The quadrature rule applies to the additional terms in the test function at time $t_{i,1}$,
which we combine with \eref{time term quadrature} to attain to the bound
\begin{multline}
        \frac{3\Delta t_i}{4}\frac12 \partial_\tau\norm{\phi(t_{i,1})}_0^2      +       \frac{\Delta t_i}{4}\frac12 \partial_\tau\norm{\phi(t_{i,2})}_0^2
                         - \frac{3\Delta t_i}{4}\frac{\Delta t_i^{-1}}{54}\norm{\Delta t_i^2\phi_{i,\tau\tau}}_0^2      \\
                \ge     \frac{\norm{\phi(t_{i^-})}_0^2  -       \norm{\phi(t_{i-1^-})}_0^2}{2}
                                + \Big(\frac{1}{72} - \frac{1}{72}\Big) \norm{\Delta t_i^2\phi_{i,\tau\tau}}_0^2        -       \hat{C}\Delta t_i \max_{0 \le j \le 2} \norm{\phi(t_{i,j})}_0^2.
\end{multline}
Accordingly, we have the bound
\begin{multline}        \label{no gronwall ineq}
        \norm{\phi(t_{i^-})}_0^2        -       \norm{\phi(t_{i-1^-})}_0^2
                        \\+     \frac{\Delta t_i}{2} \bigg\{ 3(1-\epsilon') \norm{a^{1/2}(t_{i,1})\nabla\phi_{i,TR}}_0^2                
                        + (1-\epsilon)\norm{a^{1/2}(t_{i})\nabla\phi(t_{i^-})}_0^2                                                      \\
                        +       \frac{3(1-\tilde{\epsilon}')}{32}\norm{a^{1/2}(t_{i,1}) \nabla\tilde{\phi}(\zeta_{i,1})}_0^2
                                -\frac{3(1+\tilde{\epsilon})}{32} \norm{a^{1/2}(t_{i,1})\nabla\tilde{\phi}(t_{i-1^+})}_0^2\bigg\}                                               \\
        \le     C \Delta t_i \sum_{j=1}^2 \big(        \norm{\eta_\tau(t_{i,j})}_{(-1,\mathcal{V}_h(t_{i,j}))}^2   +       \norm{\eta_{RK}(t_{i,j})}_1^2     +       \norm{\phi(t_{i,j})}_0^2  \big)
                        +       C\Delta t_i     \norm{ u(t_{i,1})  - u_{TR} (t_{i,1}) }_1^2.
\end{multline}

The negative term, $-\frac{3(1+\tilde{\epsilon})}{32} \norm{a^{1/2}(t_{i,1})\nabla\tilde{\phi}(t_{i-1^+})}_0^2$, 
prevents us from applying the discrete Gr\"onwall lemma to \eref{no gronwall ineq}.
Our approach for getting rid of this term is to create a telescoping sum once we apply the discrete Gr\"onwall lemma.
From \eref{diffusion bound} and lemma \ref{shift lemma}, and the $\mathcal{H}^1$-stability of $\mathcal{L}_2$-projection \cite{BANKYSERENTANT}, we have
\begin{multline}
        \frac{3\Delta t_i}{32} (1+\tilde{\epsilon}) \norm{a^{1/2}(t_{i,1})\nabla\tilde{\phi}(t_{i-1^+})}_0^2    \\
                \le     \frac{3\Delta t_i}{32} (1+\tilde{\epsilon})(1+\alpha \Delta t_i) \norm{\tilde{a}^{1/2}(t_{i-1})\nabla\tilde{\phi}(t_{i-1^+})}_0^2               \\
                \le     \frac{3\Delta t_i}{32} (1+\tilde{\epsilon})(1+\alpha \Delta t_i)(1+C_{\mu,d}\Delta t_i) \norm{a^{1/2}(t_{i-1})\nabla{\phi}(t_{i-1^+})}_0^2              \\
                \le     \frac{3C_{\mathcal{A},\mathcal{H}}\Delta t_i}{32}(1+\tilde{\epsilon})(1+\alpha \Delta t_i)(1+C_{\mu,d}\Delta t_i)\norm{a^{1/2}(t_{i-1})\nabla{\phi}(t_{i-1^-})}_0^2.    \\
\end{multline}
If we choose $\tilde{\epsilon} \le \frac{1}{3}$, then the controlled mesh discontinuities \eref{stability bound} and graded time stepping \eref{graded time step} gives
\(
        {3C_{\mathcal{A},\mathcal{H}}\Delta t_i}/{32}   \le     7 \Delta t_{i-1} / 8.
\)
Choose $\epsilon = 1/16$ so that
\begin{multline}
        \frac{15\Delta t_i}{16}\norm{a^{1/2}(t_{i})\nabla\phi(t_{i^-})}_0^2     -       \frac{7\Delta t_{i-1}}{8} \norm{a^{1/2}(t_{i,1})\nabla\tilde{\phi}(t_{i-1^-})}_0^2
                \\\ge   \frac{7}{8}\Big\{ \Delta t_i \norm{a^{1/2}(t_{i})\nabla\phi(t_{i^-})}_0^2 
                        -       (1+C_{\alpha,\mu,d}\Delta t_i) \Delta t_{i-1}\norm{a^{1/2}(t_{i-1})\nabla\phi(t_{{i-1}^-})}_0^2 \Big\}\\
                        +       \frac{\Delta t_i}{16} \norm{a^{1/2}(t_{i})\nabla\phi(t_{i^-})}_0^2.
\end{multline}  
Then, we use $a\ge\bar{a}$ and lemma \ref{shift lemma} to show that the bound for the $i^{\mathrm{th}}$ partition is
\begin{multline}
        \bigg[ \norm{\phi(t_{i^-})}_0^2 + \frac{7}{16}{\Delta t_i}\norm{a^{1/2}(t_{i})\nabla\phi(t_{i^-})}_0^2 \bigg]\\
                        -       (1+C'\Delta t_i) \bigg[ \norm{\phi(t_{i-1^-})}_0^2 + \frac{7}{16}{\Delta t_{i-1}}\norm{a^{1/2}(t_{i})\nabla\phi(t_{i^-})}_0^2 \bigg]    \\
                        +       \frac{\bar{a}\Delta t_i}{2} \bigg\{ 3(1-\epsilon') \norm{\phi_{i,TR}}_1^2       
                        + \frac{1}{16}\norm{\phi(t_{i^-})}_1^2   +      \frac{3(1-\tilde{\epsilon}')}{32}\norm{{\phi}(\zeta_{i,1})}_1^2 \bigg\}
        \\\le   C \Delta t_i \Big\{ \sum_{j=1}^2 \Big(         \norm{\eta_\tau(t_{i,j})}_{(-1,\mathcal{V}_h(t_{i,j}))}^2   +       \norm{\eta_{RK}(t_{i,j})}_1^2     \Big)
                        +       \norm{ u(t_{i,1})  - u_{TR} (t_{i,1}) }_1^2 +   \max_{0 \le j \le 2} \norm{\phi(t_{i,j})}_0^2   \Big\}.
\end{multline}
From the local Gr\"onwall inequality, we get the bound
\begin{multline}
        \bigg[ \norm{\phi(t_{i^-})}_0^2 + \frac{7}{16}{\Delta t_i}\norm{a^{1/2}(t_{i})\nabla\phi(t_{i^-})}_0^2 \bigg]\\
                -       (1+C'\Delta t_i) \bigg[ \norm{\phi(t_{i-1^-})}_0^2 + \frac{7}{16}{\Delta t_{i-1}}\norm{a^{1/2}(t_{i})\nabla\phi(t_{i^-})}_0^2 \bigg]    \\
                        +       \theta \Delta t_i \bigg\{ \norm{\phi_{i,TR}}_1^2                +       \norm{\phi(t_{i^-})}_1^2        +       \norm{ \phi(\zeta_{i,1}) }_1^2 \bigg\}
        \\\le   C \Delta t_i\sum_{j=1}^2 \big(         \norm{\eta_\tau(t_{i,j})}_{(-1,\mathcal{V}_h(t_{i,j}))}^2   +       \norm{\eta_{RK}(t_{i,j})}_1^2     \big)
                        +       C\Delta t_i     \norm{ u(t_{i,1})  - u_{TR} (t_{i,1}) }_1^2,
\end{multline} 
for some $\theta>0$, for $\Delta t_i$ sufficiently small.

Applying the discrete Gr\"onwall lemma gives
\begin{multline}        \label{trbdf gronwall}
        \max_{1 \le i \le m} \norm{\phi(t_{i^-})}_0^2   
                        +       \sum_{i=1}^m \Delta t_i \bigg\{ \norm{\phi_{i,TR}}_1^2          +       \norm{\phi(t_{i^-})}_1^2        +       \norm{ \phi(\zeta_{i,1}) }_1^2 \bigg\}
                \\\le   C \bigg\{       \enorm{\eta}^2 + \sum_{i=1}^m \Delta t_i \norm{u(t_{i,1})-u_{TR}(t_{i,1})}_1^2 + \norm{\phi(0)}_0^2 + \Delta t \snorm{\phi(0)}_1^2 \bigg\}.
\end{multline}
For the additional terms in the upper bound, we have
\begin{equation}
        \norm{\phi(0)}_0        \le \norm{\eta(0)}_0 \le \enorm{\eta},
\end{equation}
for $j=1,2$,
\begin{multline}
        \norm{\phi_{\tau}(t_{i,1})}_{(-1,\mathcal{V}_h(t_{i,1}))}       \\
                \le C\Big\{ \norm{\eta_{\tau}(t_{i,1})}_{(-1,\mathcal{V}_h(t_{i,1}))} + \norm{\eta_{TR}(t_{i,1})}_1 + \norm{u(t_{i,1})-u_{TR}(t_{i,1})}_1 + \norm{\phi_{TR}(t_{i,1})}_1 \Big\},
\end{multline}
\begin{equation}
        \norm{\phi_{\tau}(t_{i,2})}_{(-1,\mathcal{V}_h(t_{i,2}))}       
                \le C\Big\{ \norm{\eta_{\tau}(t_{i,2}))}_{(-1,\mathcal{V}_h(t_{i,2}))} + \norm{\eta(t_{i,2})}_1 + \norm{\phi(t_{i,2})}_1 \Big\},
\end{equation}
 and use the local Gr\"onwall lemma again to bound the maximum $\norm{\phi(t_{i,1})}_0$ at the intermediate collocation nodes to get
\begin{multline}
        \max_{\substack{1 \le i \le m\\1 \le j \le 2}} \norm{\phi(t_{i,j})}_0^2
                +       \sum_{i=1}^m \Delta t_i \bigg\{ \norm{\phi_{i,TR}}_1^2          +       \norm{\phi(t_{i^-})}_1^2        +       \norm{ \phi(\zeta_{i,1}) }_1^2
                                                \\\qquad\qquad  + \norm{\phi_\tau(t_{i,1})}_{(-1,\mathcal{V}_h(t_{i,1}))}^2     + \norm{\phi_\tau(t_{i,2})}_{(-1,\mathcal{V}_h(t_{i,2}))}^2 \bigg\}
                \\\le   C \bigg\{       \enorm{\eta}^2 + \sum_{i=1}^m \Delta t_i \norm{u(t_{i,1})-u_{TR}(t_{i,1})}_1^2 + \Delta t \snorm{\phi(0)}_1^2 \bigg\}.
\end{multline}
Since the trapezoid approximation is second order, we have
\begin{equation}
        \Delta t_i \norm{u(t_{i,1})-u_{TR}(t_{i,1})}_1^2        \le     C_{TR} \int_{t_{i-1}}^{t_i} \norm{\Delta t_i^2 u_{\tau\tau}(t)}_1^2\, dt
\end{equation}
and by the $\mathcal{H}^1\!$-stability of $\mathcal{L}_2$-projection, we have
\begin{equation}        \label{initial h1}
        \Delta t \snorm{\phi(0)}_1^2    \le C\Delta t \snorm{\eta(0)}_1^2       \le C \Delta t \Big( \snorm{\eta_{1,TR}}_1^2 + \snorm{\eta(\zeta_{1,1})}_1^2 \Big)
                \le C \enorm{\eta}^2.
\end{equation}
Thus, combining \eref{trbdf gronwall}--\eref{initial h1}, we have
\[
        \enorm{ \phi }^2       \le     C \Big\{        \enorm{\eta}^2 + \int_0^T \norm{ \Delta t_i^2 u_{\tau\tau}(t) }_1^2\, dt \Big\},
\]
as desired.
\end{proof}

\section{Numerical Experiments}\label{sec4}

In this section, we present some  numerical examples that illustrate TR-BDF2  quadratic moving finite elements.
Although Theorem \ref{trbdf a priori theorem} assumed for simplicity the collocation node $\hat{t}_1  =2/3$, 
corresponding to Gauss-Radau quadrature, careful examination of our proof, as well as previous experiments \cite{METTITHESIS}
verify that second order convergence in time holds for nearby values of this collocation node.
The best performance we observed was for the special value $\hat{t}_1= 2-\sqrt2$, as anticipated
\cite{A30,A30a,choudhury1992waveform,ying2009composite}.
Accordingly, the PDEs are solved using this collocation node.

A solver for time-dependent linear convection-diffusion-reaction equations with a single dimensional space domain was written in C++.
For simplicity, an approximate method of characteristics was used to generate the mesh motion;
we used two steps of forward Euler to generate quadratic trajectories that approximately satisfy $x_t = b(x)$ at each spatial node.
Between the time partitions, we reconfigure the mesh to be a uniform partition of the spatial domain to avoid mesh degeneration.

If any mesh nodes collide or spread too far apart, one of the nodes is deleted or the element is bisected, respectively.
We use interpolation rather than $\mathcal{L}_2$-projection at these mesh discontinuities to reduce CPU time.
The effects of this modification are briefly discussed below and in \cite{METTITHESIS}.
Uniform time steps are used and the time domain remains fixed in these experiments so that more time steps implies smaller $\Delta t$ 
rather than a longer simulation.
All experiments in this section report the accuracy and CPU time of the solution computed on the moving 
mesh relative to a solution on a uniform and non-moving mesh.

We first test our methods on two linear convection-diffusion problems, both of which live on the domain $(x,t) \in (-3,3)\times(0,1]$.
The first problem is a convection-dominated problem given by
\begin{equation}        \label{convection dominated experiment}
        u_{1,t}(x,t) - 0.01 u_{1,xx}(x,t) + 3 u_{1,x}(x,t)      = f_1(x,t),
\end{equation}
with $f_1$, the initial condition, and Neumann boundary condition chosen such that the solution is given by
\[
        u_1(x,t)        = e^{-(x-3t)^2}.
\]
Each time partition is initialized with a uniform mesh and since the convection velocity is constant in this problem,
the characteristic trajectories perfectly cancel out its effect, $x_t = 3$.
A standard Galerkin {approximation} of the differential equation is imposed for both the moving and non-moving meshes at the discrete collocation nodes,
as in \eref{midstep constraint}--\eref{endstep constraint}.

\begin{figure}
        \centerline{
        \includegraphics[scale=.5]{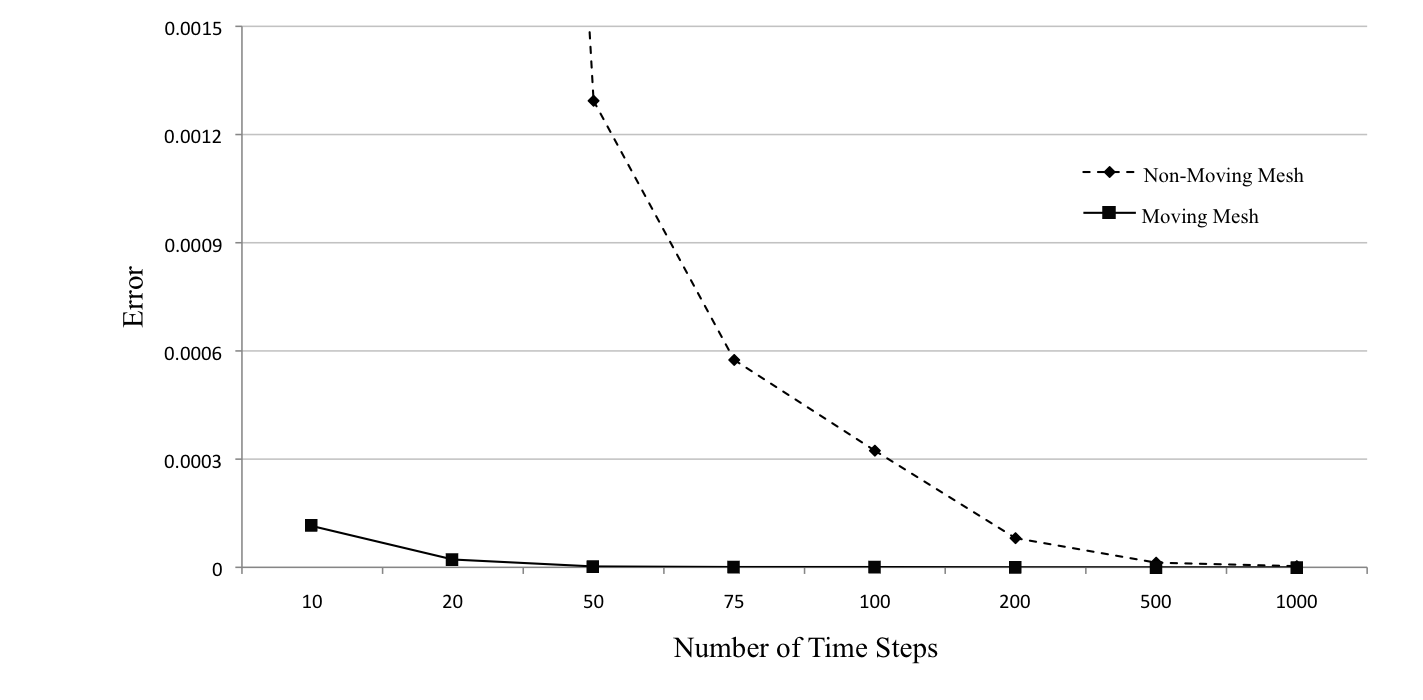}
        }
        \caption{The final $\mathcal{L}_2$-error at time $t=1$ of solutions computed for \eref{convection dominated experiment} on non-moving (dashed line) and moving (solid line) meshes.
        The error is significantly reduced using a moving mesh, particularly for larger time steps.}
        \label{convection error plot}
\end{figure}

Problem \eref{convection dominated experiment} sees a great advantage when the mesh moves with the convection velocity,
with up to a 300-fold improvement in the relative error on highly refined meshes with large time steps.
The ratio of the final $\mathcal{L}_2$-errors of the moving mesh and static mesh solutions at the end of the simulation, 
\[
        \mathcal{L}_2\mathit{-error} = \norm{u(\cdot,1)-\tilde{u}(\cdot,1)}_0/ \norm{u(\cdot,1)-\bar{u}(\cdot,1)}_0,
\]
 and CPU times are reported in Table \ref{convection basic motion table},
 where $\tilde{u}$ and $\bar{u}$ respectively represent the moving and static mesh solutions.
These numbers are relative to the solution computed using a static mesh;
values less than 1 correspond to an decrease in the norm of the error or a speedup in the CPU time, when using moving meshes.
The accuracy of the solution at the end of the simulation is plotted for $n=1001$ with respect to the number of time steps in figure \ref{convection error plot}.

The greatest gains are realized when $\Delta t$ is large relative to $\Delta x$, validating moving finite elements improved stability for larger time steps.
However, notice that the use of coarse meshes in space coupled with short time steps leads to unsatisfactory performance.
This is caused by large interpolation errors accumulating over many mesh discontinuities.
Experiments in \cite{METTITHESIS} demonstrate that these instabilities can be suppressed by using $\mathcal{L}_2$-projection,
as required by the analysis of Section \ref{sec3}, for such 
discretizations---since the mesh is coarse in these cases, the speedup from interpolation is rather subtle, which justifies using projection.

\begin{table}
        \label{convection basic motion table}
        \begin{tabular}{| l | c | c | c | c | c | c | c | c | c |}
                \hline
                        &       \multicolumn{2}{| c } {$n=101$}          & \multicolumn{2}{| c } {$n=501$}
                        &       \multicolumn{2}{| c } {$n=1001$}         & \multicolumn{2}{| c |} {$n=3001$}    \\
                \hline
                $m$     &       $\mathcal{L}_2$-error   &       CPU             &       $\mathcal{L}_2$-error   &       CPU
                        &       $\mathcal{L}_2$-error   &       CPU             &       $\mathcal{L}_2$-error   &       CPU     \\
                \hline
        10              &0.0096 &1.018  &0.0049 &1.009  &0.0036 &1.007  &0.0031 &1.014  \\
        20              &0.0116 &1.047  &0.0031 &1.044  &0.0027 &1.041  &0.0018 &1.048  \\
        50              &0.0235 &1.066  &0.0031 &1.062  &0.0019 &1.060  &0.0016 &1.072  \\
        75              &0.2559 &1.067  &0.0031 &1.059  &0.0016 &0.576  &0.0031 &1.063  \\
        100             &0.0471 &1.072  &0.0018 &1.064  &0.0029 &0.521  &0.0042 &1.070  \\
        200             &4.8120 &1.068  &0.0304 &1.070  &0.0036 &1.062  &0.0065 &1.064  \\
        500             &23.420 &1.073  &0.0037 &1.068  &0.0048 &1.063  &0.0083 &1.070  \\
        1000            &30.122 &1.079  &3.8142 &1.079  &0.0047 &1.079  &0.0087 &1.074  \\
                \hline
        \end{tabular}
        \caption{A comparison table for problem \eref{convection dominated experiment} between non-moving and moving mesh methods using interpolation.
        Reported for each mesh discretization is the ratio of the $\mathcal{L}_2$-errors at time $t=1$
        and the relative increase in CPU time of the moving mesh solution to the static mesh solution.
        Moving the mesh leads to clear advantages when the spatial mesh is sufficiently refined,
        though buildup in error occurs when the time discretization is too refined, relative to the spatial discretization.}
\end{table}

\begin{figure}
        \centerline{
        \includegraphics[scale=.5]{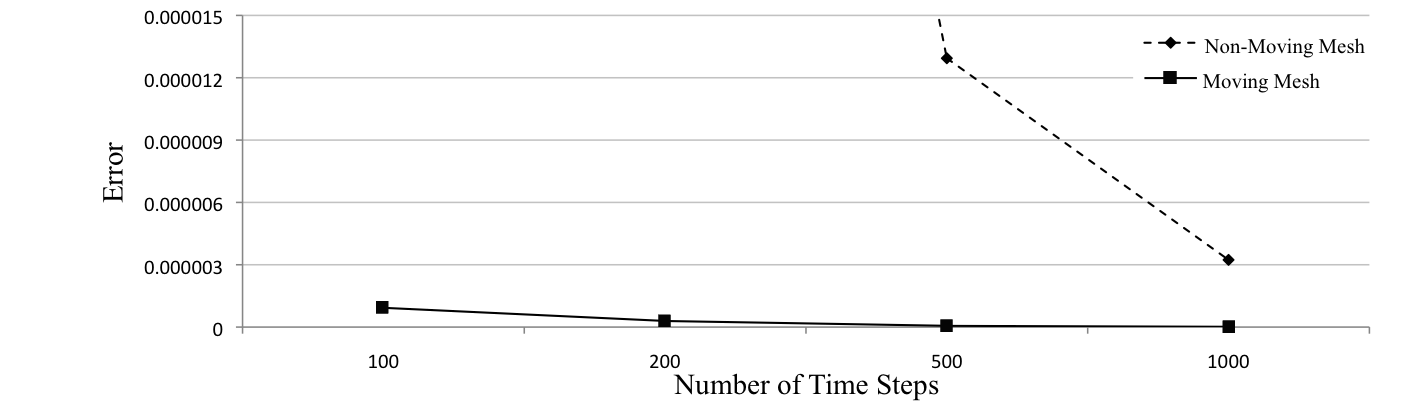}
        }
        \caption{A close-up of figure \ref{convection error plot}.  Even for highly refined partitions of the time domain, the moving mesh solution strongly outperforms the static mesh solution.}
        \label{convection error plot zoom}
\end{figure}

The next problem we consider is given by
\begin{equation}        \label{diffusion experiment}
        u_{2,t}(x,t) - \big( (x^2 + t^2 + 0.1)\,u_{2,x}(x,t)\big)_x + 0.1(x^3-9x) u_{2,x}(x,t) + u_2(x,t)       = f_2(x,t),
\end{equation}
and the source term and boundary condition are chosen so that the solution of the differential equation is given by
\[
        u_2(x,t)        = \sin \Big( \frac{\pi}{6}(x+5t) \Big).
\]
Notice that the magnitude of the diffusion is comparable to the convection velocity, so the effects of the moving elements is primarily to counterbalance
the asymmetry of the PDE brought on by the convection term.

A slight advantage is recognized by using moving meshes, when $\Delta t$ is not too small relative to $\Delta x$.
A comparison of the errors and CPU time is given in Table \ref{diffusion basic motion table}.
Comparing the results of this experiment to those of problem \eref{convection dominated experiment},
it is clear that this problem does not benefit nearly as much from a moving mesh nor is the accuracy as dramatically impacted when the time discretization is
much finer than the spatial discretization.
Since the diffusion term had a stronger presence in this problem, the solution computed on a non-moving mesh demonstrates comparable accuracy,
though small gains come from moving the mesh due to the cancellation of the mesh and convection velocities.

\begin{table}
        \label{diffusion basic motion table}
        \begin{tabular}{| l | c | c | c | c | c | c | c | c | c | c |}
                \hline
                        &       \multicolumn{2}{| c } {$n=101$}  & \multicolumn{2}{| c } {$n=501$}
                        &       \multicolumn{2}{| c } {$n=1001$}        & \multicolumn{2}{| c |} {$n=3001$}     \\
                \hline
                $m$     &       $\mathcal{L}_2$-error   &       CPU     &       $\mathcal{L}_2$-error   &       CPU 
                        &       $\mathcal{L}_2$-error   &       CPU     &       $\mathcal{L}_2$-error   &       CPU             \\
                \hline
        10              &0.8012 &1.072  &0.8014 &1.073  &0.8014 &1.082  &0.8014 &1.070  \\
        20              &0.7927 &1.068  &0.7928 &1.075  &0.7928 &1.053  &0.7928 &1.071  \\
        50              &0.7912 &1.070  &0.7893 &1.072  &0.7893 &1.068  &0.7893 &1.072  \\
        75              &0.7946 &1.069  &0.7888 &1.072  &0.7888 &1.070  &0.7888 &1.084  \\
        100             &0.7997 &1.071  &0.7886 &1.072  &0.7886 &1.060  &0.7886 &1.066  \\
        200             &0.8298 &1.070  &0.7895 &1.082  &0.7884 &1.071  &0.7885 &1.066  \\
        500             &0.9405 &1.072  &0.8009 &1.071  &0.7908 &1.071  &0.7883 &1.067  \\
        1000            &0.9791 &1.063  &0.8387 &1.069  &0.8010 &1.072  &0.7884 &1.072  \\
                \hline
        \end{tabular}
        \caption{A comparison table for problem \eref{diffusion experiment} between static mesh methods and moving mesh methods using interpolation.
        Reported for each mesh discretization is the ratio of final $\mathcal{L}_2$-errors and the relative increase in CPU time of the moving mesh solution to the static mesh solution.}
\end{table}

\begin{figure}
        \centerline{
        \includegraphics[scale=.5]{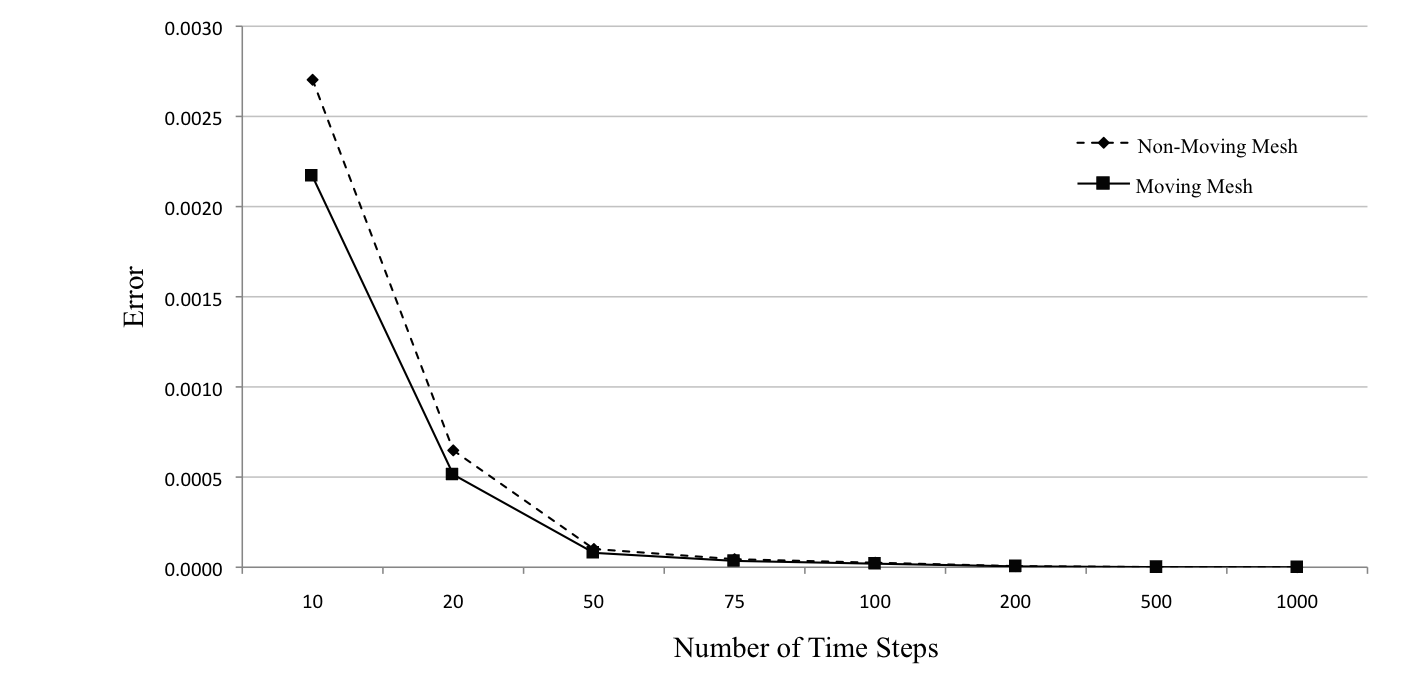}
        }
        \caption{The final $\mathcal{L}_2$-error at time $t=1$ of solutions computed for problem \eref{diffusion experiment} on non-moving (dashed line) and moving (solid line) meshes.
        Small gains are made when the mesh moves, especially there are few time steps, as $\Delta t$ is large.}
        \label{diffusion error plot}
\end{figure}

In problem \eref{diffusion experiment}, setting $x_t=5$ would move the mesh along the characteristic trajectories of the solution,
rather than canceling the convection velocity.
This suggests that there are other mesh motion schemes that may improve the accuracy of the computed solution other than the method of characteristics.
Other schemes have been proposed that compute the mesh motion and solution as a coupled system of (potentially nonlinear) ODEs,
though many of these do not fit into the framework of the moving finite element scheme in Section \ref{sec3} due to their nonlinearity
\cite{CARLSONMILLER1,CARLSONMILLER2}.
Other schemes use error estimates \cite{adjerid1986moving1d,adjerid1986moving}, predictor-corrector techniques \cite{baines1994moving},
 and conservation laws to determine the mesh motion \cite{baines2005moving,baines2011velocity}.
A mesh motion scheme of primary interest could use a posteriori error estimates and adaptive meshing
($h$-refinement at the mesh discontinuities and mesh smoothing to evolve the mesh continuously) to find an appropriate mesh.
Such a solver can be implemented leveraging the adaptive meshing routines of some existing software like PLTMG \cite{bank2012pltmg}.

One final experiment that we consider is designed to evaluate the moving finite element scheme for a nonlinear problem.
We apply our moving space-time method to Burgers' equation with a large Reynolds number.
Burgers' equation is an interesting test problem for our scheme as it is a simple nonlinear equation that develops steep moving fronts that sweep through the domain.
Artificial oscillations are commonly found in the computed solution near the shock layer when the time discretization is not sufficiently refined \cite{MILLER1,MILLER2}.
The differential equation given by
\begin{equation}        \label{burgers equation}
        u_t     -       \frac{1}{R} u_{xx}      +       uu_x            = 0     \quad\mathrm{for}\ (x,t)\ \mathrm{in}\ (-3,3)\times(0,2],
\end{equation}
where we choose a large $R>0$, and we assume Neumann boundary conditions, $u_x(\pm3)=0$.
The initial condition is chosen so that a moving front forms in the middle of the domain and propagates toward the right boundary.

We use the method of characteristics, which gives $x_t \approx u$.
The solution is computed at the end of each time step so that we can use linear mesh motion $x_t(t) =  u_h(x,t_{i-1^-})$ for $t_{i-1} \le t \le t_i$.
Furthermore, we use a single step of Newton's method to solve the nonlinear equation at each collocation node.
Unlike the earlier experiments, 
we do not reset the mesh to be uniform at the beginning of each time partition.
This allows the nodes to appropriately accumulate near the shock layer, where nodes are deleted once they get too close together.
In addition to the standard Galerkin discretizations for the non-moving and moving meshes,
we also compare our results to a solution computed using a Streamline Upwind/Petrov Galerkin (SUPG) discretization on a fixed mesh,
as in \cite{SUPG}.

\begin{figure}
        \centerline{
        \includegraphics[scale=.6]{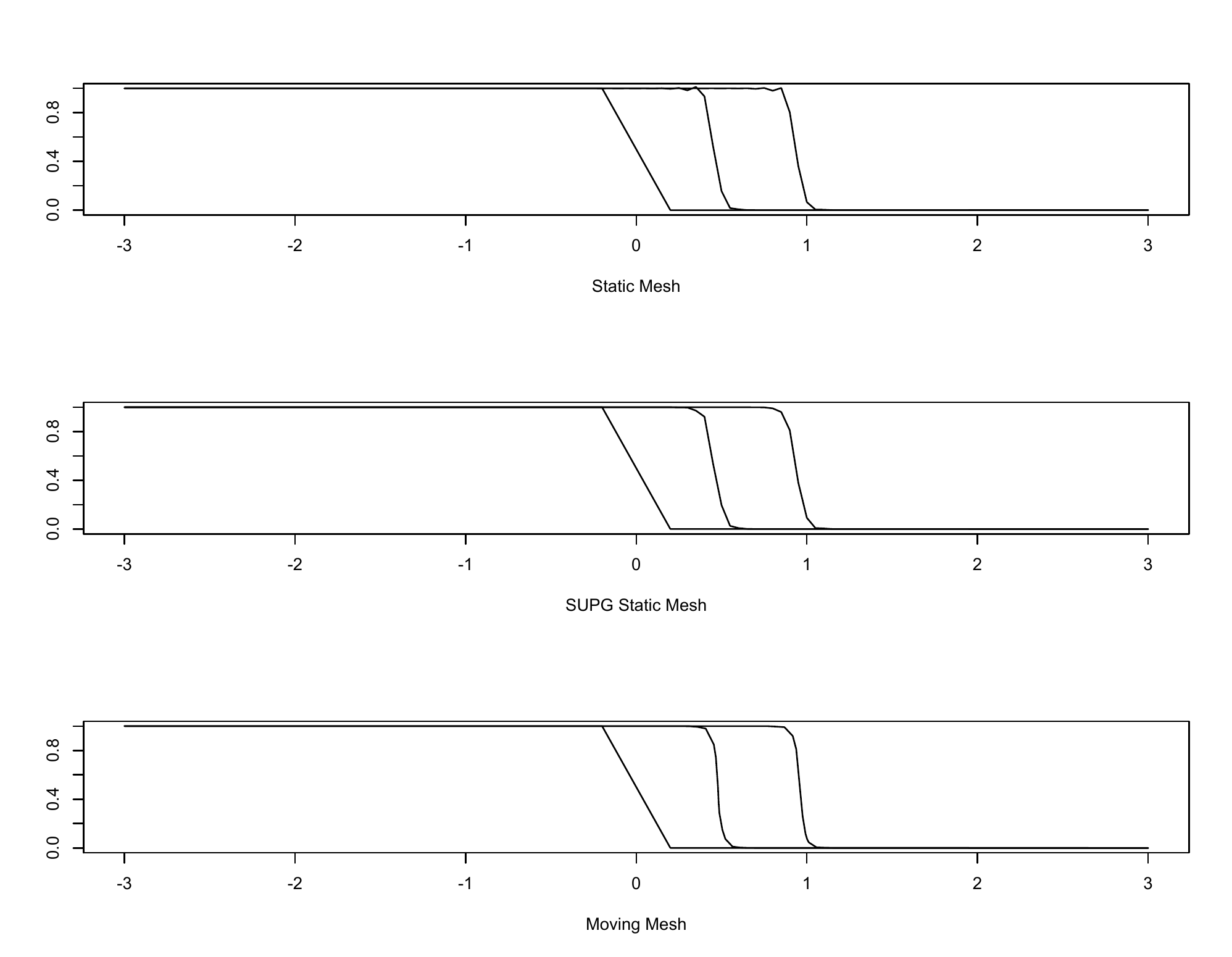}
        }
        \caption{Solutions computed for Burgers' equation with $R=100$, using $n=61$ spatial nodes and $m=25$ time steps.
        The top graph depicts the solution computed on a non-moving mesh, where oscillations develop near the shock layer.
        The middle graph plots the solution computed using SUPG, using a SUPG coefficient of $0.1$.
        The bottom graph shows that the moving finite element solution has successfully dampened these spurious oscillations.
        In comparing the moving finite element solution with the solution computed using SUPG,
        the moving finite element solution maintains a steeper shock layer, whereas artificial diffusion has smeared the layer in the SUPG solution.}
        \label{thick shock fig}
\end{figure}

\begin{figure}
        \centerline{
        \includegraphics[scale=.5]{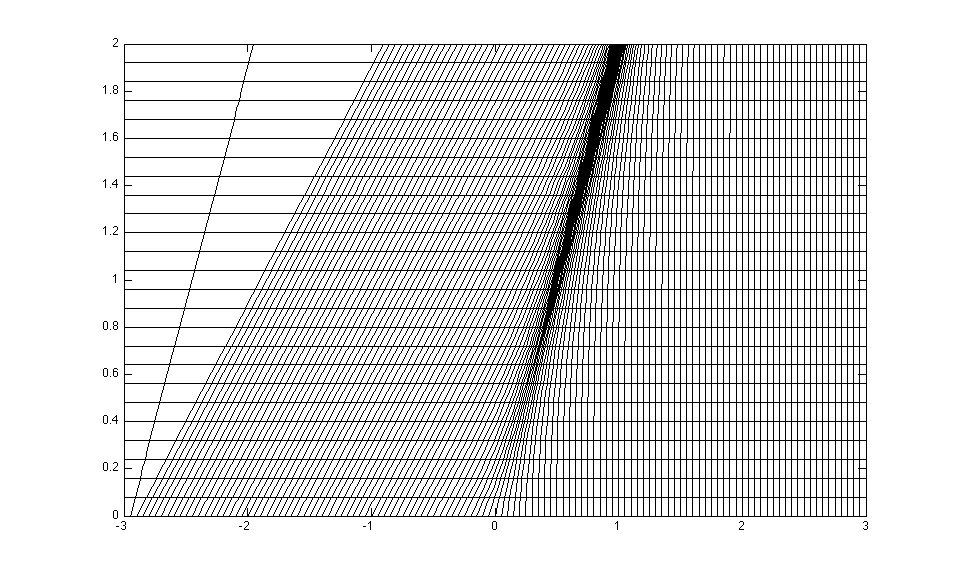}
        }
        \caption{An example of a moving mesh with $m=25$ time steps and initialized with $n=61$ spatial nodes at the beginning of the simulation.
        The method of characteristics sets $x_t \approx u(x)$, where both the hat and bump node trajectories are depicted.
        Comparing the density of the spatial nodes to the figures of the computed solutions in figure \ref{thick shock fig},
        the spatial nodes properly congregate near the shock layer and track it throughout the simulation.}
        \label{burger mesh figure}
\end{figure}

\begin{figure}
        \centerline{
        \includegraphics[scale=.6]{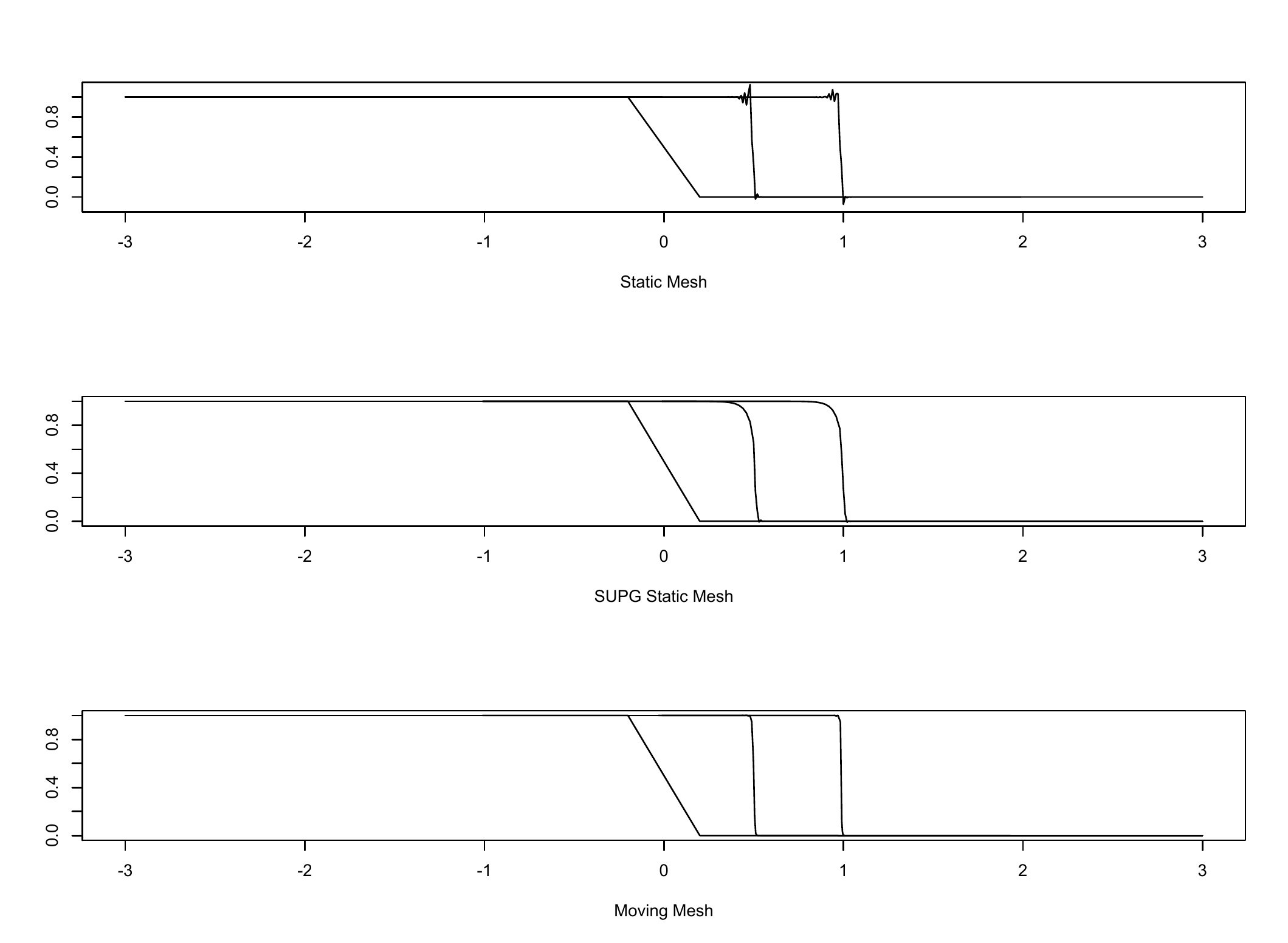}
        }
        \caption{Solutions computed for Burgers' equation with $R=1000$, using $n=301$ spatial nodes and $m=100$ time steps.
        The top graph depicts the solution computed on a non-moving mesh, where large oscillations develop near the shock layer.
        The middle graph shows the solution computed using SUPG, where the SUPG coefficient is set to $1$.
        The top of the shock layer in the solution computed using SUPG is not sharply defined and small oscillations are present at its foot.
        The bottom graph shows that the moving finite element solution has almost completely dampened these spurious oscillations,
        while preserving the sharpness of the solution on both sides of the shock layer.}
        \label{thin shock fig}
\end{figure}

Figure \ref{thick shock fig} displays solutions for equation \eref{burgers equation} with $R=100$,
computed on fixed meshes with standard Galerkin and a SUPG discretizations in space,
where the SUPG coefficient of $0.1$, as well as a moving finite element discretization.
Numerically-induced oscillations are expectedly present in solution found on the non-moving mesh \cite{XUHU},
whereas they have already been suppressed in the SUPG and the moving finite element solutions.
From figure \ref{thick shock fig}, the moving finite element solution maintains a steep drop into the moving front,
where artificial diffusion is clearly present at the top of the shock layer in the SUPG solution.
This demonstrates the increased flexibility in the time discretization when moving the mesh,
without the artificial diffusion of a SUPG discretization.
Experiments combining SUPG and moving meshes give similar results to the standard moving finite element discretization
since the ``streamline direction,'' now given by $u-x_t$, is small.
Figure \ref{burger mesh figure} depicts an example of a moving mesh used for solving the PDE,
where we see the spatial nodes accumulating at the shock layer, as desired.
We also ran simulations where the Reynolds number is set to $R=1000$.
This reduces the diffusive forces in the equation and leads to a thinner shock.
Solutions computed with 100 time steps are displayed in figure \ref{thin shock fig},
where we again see the moving finite element has remarkably smaller oscillations near the shock layer without the
undesirable effects of artificial diffusion.

\section{Conclusion} \label{sec5}

Theorem 4.3 of \cite{BANKMETTI} provides a symmetric error bound for some space-time moving finite element methods that use a one parameter family of collocation methods.
Employing TR-BDF2 for time integration, however, does not fit into this framework and complicates the analysis of this finite element method, 
ultimately breaking the symmetry of the error bound.
Nevertheless, it has been shown that Theorem \ref{trbdf a priori theorem}  analytically implies second order accuracy with respect to mesh refinement \cite{METTITHESIS}
and the moving finite element scheme yields superior performance over the method of lines approach with a stationary mesh,
as demonstrated by the numerical results.
Furthermore, improved accuracy for the moving finite element solution was maintained in application to a simple nonlinear problem,
although the theoretical analysis does not cover such cases.

In the numerical experiments, the method of characteristics was used to determine the mesh motion;
however, an optimally robust and well-defined mechanism for evolving the mesh for general PDEs is still a matter of active research.
It is suspected that predictor-corrector methods, coupled with adaptive meshing for the spatial discretization,
can be powerful tools in moving the mesh without requiring additional user-supplied information about a PDE or its solution.

\bibliography{ref}
\bibliographystyle{siam}
\end{document}